\RequirePackage{fix-cm}
\documentclass[smallextended]{svjour3}       
\smartqed  
\usepackage{graphicx}
\usepackage{mathptmx}      
\usepackage{amsmath}
\usepackage{amssymb}
\usepackage{algorithm}
\usepackage[noend]{algpseudocode}
\usepackage[hidelinks,pdfusetitle,pdfdisplaydoctitle]{hyperref}
\usepackage{booktabs}
\usepackage{xcolor}
\usepackage[inline]{enumitem}

\allowdisplaybreaks
\newcommand{\Set}[1]{\mathcal{#1}}
\newcommand{\R}{\mathbb{R}}
\newcommand{\Z}{\mathbb{Z}}

\journalname{}
\begin{document}

\title{A Decision Diagram Approach for the Parallel Machine Scheduling Problem with Chance Constraints}

\titlerunning{A DD approach for the CC-PMSP}        

\author{Nicolás Casassus \and Margarita Castro \and Gustavo Angulo}


\institute{N. Casassus \at
              Department of Industrial and Systems Engineering, Pontificia Universidad Católica de Chile
           \and
           M. Castro \at
              Department of Industrial and Systems Engineering, Pontificia Universidad Católica de Chile \\
              \email{margarita.castro@uc.cl}           
           \and
           G. Angulo \at
              Department of Industrial and Systems Engineering, Pontificia Universidad Católica de Chile
}

\date{Received: date / Accepted: date}

\maketitle
\begin{abstract}
The Chance-Constrained Parallel Machine Scheduling Problem (CC-PMSP) assigns jobs with uncertain processing times to machines, ensuring that each machine's availability constraints are met with a certain probability. We present a decomposition approach where the master problem assigns jobs to machines, and the subproblems schedule the jobs on each machine while verifying the solution's feasibility under the chance constraint. We propose two different Decision Diagram (DD) formulations to solve the subproblems and generate cuts. The first formulation employs DDs with a linear cost function, while the second uses a non-linear cost function to reduce the diagram's size. We show how to generate no-good and irreducible infeasible subsystem (IIS) cuts based on our DDs. Additionally, we extend the cuts proposed by Lozano \& Smith \cite{lozano2018} to solve two-stage stochastic programming models. Our DD-based methodology outperforms traditional integer programming (IP) models designed to solve the CC-PMSP in several instances. Specifically, our best DD-based approach solves 55 more instances than the best IP alternative (from a total of 405) and typically achieves smaller gaps (50\% vs. 120\% gap on average). 
\end{abstract}

\section{INTRODUCTION}
\label{sec:Introduction}
Scheduling problems represent a fundamental area within operational research, with numerous practical applications across economic and industrial sectors \cite{pinedo2023}. These problems focus on efficiently allocating limited resources over time to execute a set of tasks while optimizing one or more objectives. They naturally arise in various contexts, such as manufacturing environments where jobs must be assigned to machines and their processing order determined, hospitals scheduling surgeries and medical interventions, airports establishing takeoff and landing sequences, and logistics centers coordinating loading, unloading, and distribution processes, to name a few \cite{baker2019,blazewicz2019,brucker2006,cardoen2010,rinnooy1976,lenstra1977,pinedo2023}. Efficient scheduling of these activities can significantly contribute to cost reduction, optimize the use of limited resources, reduce waiting times, and improve service quality.

Among the various scheduling problems, the Parallel Machine Scheduling Problem (PMSP) involves scheduling a set of jobs on a set of machines, considering that each job has an execution time and a setup time between jobs \cite{baker2019,brucker2006,pinedo2023}. Multiple objective functions can be used to optimize the problem, such as minimizing job execution time, minimizing job tardiness, and minimizing setup times. This problem, classified as NP-hard, has been extensively studied in scientific literature due to its combinatorial complexity and practical relevance \cite{baker2019,brucker2006,lenstra1977}. However, traditional deterministic models assume perfect knowledge of processing and setup times, an assumption that is rarely met in real-life applications. In practice, execution times are subject to variations caused by multiple factors, including equipment breakdowns, changes in resource availability, and variability in human performance, among others. This uncertainty can significantly compromise the feasibility and quality of solutions obtained through deterministic approaches.

This research focuses on the Chance-Constrained Parallel Machine Scheduling Problem (CC-PMSP), a version of PMSP that incorporates uncertainty in job execution and setup times. The problem seeks to assign jobs to different machines, guaranteeing with a certain probability that all assigned jobs can be completed within the available machine time. Given the computational challenges of solving the CC-PMSP, we propose a decomposition approach that divides the problem into two levels. The master problem handles job assignment to machines, determining which jobs will be processed by each machine. The subproblems, one for each machine and scenario, determine the optimal processing sequence of the assigned jobs and verify the feasibility of the solution under the chance constraint.

One of the main drawbacks of such decomposition is that the subproblems correspond to NP-hard problems (i.e., sequence problems with setup times). This work proposes Decision Diagrams (DDs), which are graphical structures employed in optimization to represent the feasibility set of discrete optimization problems \cite{bergman2016decision}, to solve these subproblems, aiming to leverage their high reusability across different subproblems. These DDs have been utilized in numerous applications, particularly in scheduling and sequencing problems, with high success \cite{cire2013}. Our decomposition encodes the subproblem using a DD and generates cuts from it to guide the master problem toward optimal solutions. We propose two different DD formulations to address this challenge. The first one, called DD-LastJob (DD-LJ), employs a linear cost function to evaluate paths in the diagram. The second, called DD-JobSet (DD-JS), employs a non-linear cost function that significantly reduces the diagram's size. Additionally, we develop mechanisms to generate No-Good and Irreducible Infeasible Set (IIS) cuts based on our DDs, which efficiently eliminate infeasible solutions from the search space.

In summary, the main contributions of this work are the following:

\begin{enumerate}
    \item We present a novel variant of the CC-PMSP and propose a decomposition methodology based on DDs.
    \item We propose two DD formulations (DD-LJ and DD-JS) to efficiently solve the subproblems and present novel algorithms to generate  No-Good and IIS cuts based on such DDs.
    \item We adapt the methodology proposed by Lozano \& Smith \cite{lozano2018} for two-stage stochastic programming problems to the CC-PMSP, which can be easily extended to other chance-constrained problems.
    \item We provide extensive numerical experiments and empirical evidence demonstrating that the proposed DD-based methodology outperforms Integer Programming (IP) alternatives.
\end{enumerate}

The remainder of the paper is organized as follows.  Section \ref{sec:LiteratureReview} establishes our research foundation with a literature review, followed by the description and formalization of the CC-PMSP in Section \ref{sec:DescriptionFormulation}. Section \ref{sec:decomposition} introduces our decomposition strategy, while Section \ref{sec:DDModels} reviews DDs and presents the DD-LJ and DD-JS formulations for CC-PMSP subproblems. Section \ref{sec:DDCuts} details algorithms for generating No-Good and IIS cuts based on these DDs, followed by our computational experiments and results analysis in Section \ref{sec:Experiments}. Finally, Section \ref{sec:Conclusions} presents our findings, draws conclusions, and outlines directions for future research.

\section{Literature Review}
\label{sec:LiteratureReview}

This research studies chance-constrained programming (CCP) problems with binary variables. Specifically, we consider  problems with the following structure
\begin{subequations}
    \begin{alignat}{2}
        \max && \mathbf{u}^{\top} \mathbf{x} \label{FO:1}\\
        \text{s.t.} && \mathbf{A}\mathbf{x} &\leq \mathbf{b} \label{Constr:1} \\
        && \mathbb{P}(\exists \mathbf{y} \in \Set{Y}(\mathbf{x}, \omega)) &\geq 1 - \varepsilon  \label{ChanceConstr:1}\\
        && \mathbf{x} &\in \{0,1\}^n , \label{Binary:1} 
    \end{alignat}
\end{subequations}
where the objective function \eqref{FO:1} seeks to maximize the sum of utilities associated with the solution, subject to deterministic linear constraints \eqref{Constr:1}, a chance constraint \eqref{ChanceConstr:1} ensuring feasibility with probability at least $1 - \varepsilon$ ($\varepsilon\in (0,1)$) with respect to realizations $\omega\in\Omega$ of random parameters, and binary constraints \eqref{Binary:1}. We now review the literature related to CCP problems and DDs.

\subsection{Chance-Constrained Programming Problems}
\label{subsec:ChanceConstrLiteratureReview}

CCP problems have been studied since at least the work of \cite{prekopa1970}, but remain computationally challenging for two main reasons: (i) the solution space is not convex, and (ii) evaluating the feasibility of a solution requires integrating multi-dimensional functions \cite{luedtke2014}. Traditional methods only achieve good solutions in special cases. For example, if the chance constraint is a linear inequality where the random coefficients follow a normal distribution, it is possible to reformulate the problem as a deterministic non-linear convex problem if $\varepsilon < 0.5$ \cite{calafiore2006}. Another example is for problems with $\Set{Y}(\mathbf{x}, \omega) = \{\mathbf{y}: T\mathbf{x} + W\mathbf{y} \geq h(\omega)\}$, where $h(\omega)$ follows a discrete finite distribution, which can be solved with a IP formulation \cite{kucukyavuz2012}.

Luedtke \cite{luedtke2014} proposed a formulation to exactly solve problems with chance constraints, considering discrete distributions with finite support and random polyhedral constraints. Their procedure decomposes the problem into a master problem containing all the deterministic constraints and one subproblem for each scenario, and then solves the problem using a branch-and-cut algorithm. Although general, this procedure is typically tailored to each application, depending on the structure of the subproblem and the nature of its variables. Liu et al. \cite{liu2016} extended this formulation to include second-stage decisions with additional costs and consideration of the magnitude of not meeting certain scenarios. Building on these foundations, Elçi \& Noyan \cite{elci2018} further extended these formulations to optimize humanitarian aid distribution and capacity planning for future disasters by incorporating risk measures through CVaR. Similarly, Lodi et al. \cite{lodi2019} adapted these approaches by relaxing problem assumptions to enable application in hydroelectric generation planning.

To the best of our knowledge, our work is the first to consider chance constraints in the PSMP with uncertainty in the processing and setup times. Song \& Leus \cite{song2022parallel} defined a variant of the CC-PMSP with uncertain processing times following a normal distribution, and does not consider setup times. Thus, their chance constraints can be reformulated as a second-order conic constraint, which can be handled by commercial solvers. The authors proposed a branch-and-cut approach and dichotomic search procedures to solve the problem more efficiently. In contrast, our CC-PMSP variant considers uncertain processing and setup times drawn from any random distribution with finite support, which makes the problem considerably more challenging, as it incorporates a joint-chance constraint for each machine to address the sequential aspect of the problem. 

Researchers have also applied chance constraints to handle uncertainty in other related scheduling problems. Deng \& Sheng \cite{deng2016decomposition} address a multi-server appointment scheduling problem with chance constraints where appointment times are uncertain with finite support. Similarly to our work, they proposed a two-stage decomposition to efficiently solve the problem, where their decomposition strongly relies on the packing structure of their problem. Additionally, similar to \cite{song2022parallel}, the authors do not consider uncertain setup times, which simplifies their subproblems compared to ours. Additional related works include applications of chance constraints to scheduling problems, such as project scheduling  \cite{wang2017}, surgery planning  \cite{noorizadegan2018}, workforce scheduling \cite{excoffier2015}, job-shop scheduling \cite{shen2016}, and multi-depto vehcile scheduling \cite{castro2024incorporating}, to name a few. While these applications are somewhat related to our problem, we note that none of these works consider a joint chance constraint similar to the one proposed for our CC-PMSP variant. 


\subsection{Decision Diagrams}
\label{subsec:DDLiteratureReview}

Decision Diagrams (DD) are graph-based structures with many applications. Initially introduced for circuit design \cite{hu1996,lee1959}, DDs have also proven helpful for sequential pattern mining and genetic programming \cite{loekito2010,wegener2000}. In optimization, DDs have been shown to be useful and versatile tools for solving various discrete problems with mathematical programming techniques \cite{bergman2016,castro2022}. For instances, DDs can be employed in cut generation algorithms \cite{becker2005,behle2007b,castro2021,davarnia2021,tjandraatmadja2019}, in global constraints \cite{hooker2006}, in bound computations \cite{bergman2014,cappart2019,castro2020,vanhoeve2022,maschler2021,rudich2022} and in handling stochastic problems \cite{latour2019,lozano2018}.

Several studies utilize DDs to solve scheduling problems, beginning with Cire \& van Hoeve \cite{cire2013} who proposed a general approach for the sequencing problem. Since this seminal work, several researchers have employed DDs to solve scheduling and sequence problems (e.g., vehicle routing variants \cite{kinable2017hybrid,castro2020mdd}, and clinical rotation optimization \cite{cire2019}, to name a few). Indeed, several works have tackled deterministic variants of the PMSP. Lozano et al. \cite{lozano2020} to solve PMSP, where jobs consist of two tasks with a minimum time interval between them. Also, Kowalczyk \& Leus \cite{kowalczyk2018} implemented branch-and-price based on DDs for scheduling problems with parallel machines, and Raghunathan et al. \cite{raghunathan2022} used a similar algorithm to optimize time and number of transfers for first and last mile of a transportation system. We refer the reader to the survey \cite{castro2022} and tutorial \cite{vanhoeve2024} for more DD applications.

While DDs have been extensively used in deterministic problems, there is limited work on solving problems with uncertain parameters. Lozano \& Smith \cite{lozano2018} proposed a solution approach to a special type of two-stage stochastic integer problems using DDs with their novel incorporation of capacitated arcs. This idea enables parameterizing arc availability as a function dependent on a binary vector $\mathbf{x}$ and generates Benders' cuts in a standard two-stage decomposition approach. MacNeil \& Bodur \cite{macneil2024} extended this DD-based approach, allowing for a greater variety of two-stage stochastic problems where the first stage influences second-stage constraints. They argue that using capacitated arcs can lead to poor DD reductions and propose a new parameterization based on arc costs. Guo et al. \cite{guo2021} adopted the cuts presented by \cite{lozano2018} to propose a decomposition to solve the Stochastic Distributed Operating Room Scheduling problem.

To the best of our knowledge, there are only two other works that have utilized DD for CCP problems. Latour et al. \cite{latour2017,latour2019} employ DDs to represent  decisions and the probabilistic distributions of the problem. Their DD is reformulated into a quadratic constraint model, which is then linearized and introduced into an IP formulation of the original problem. Haus et al. \cite{haus2017} also use DDs to solve two-stage stochastic problems where first-stage decisions affect the stochastic process (i.e., endogenous uncertainty), which are then used to reformulate the chance constraints into linear inequalities. To our knowledge, no other techniques have been developed to solve scheduling problems with chance constraints similar to the CC-PMSP proposed in this work (i.e., a decomposition approach where subproblems are solved using DDs).

\section{Problem Description and Formulation}
\label{sec:DescriptionFormulation}

The CC-PMSP involves deciding the assignment and order of a set $\Set{J}$ of jobs on a set $\Set{M}$ of machines. Each job $j \in \Set{J}$ has an assigned utility $f_j$ that dictates its importance in scheduling, as well as execution times $t_j$ and setup times $d_{jk}$ with respect to another job $k \in \Set{J}$. These times, both execution and setup, are stochastic and are not known at the time of optimization. We will identify this stochasticity with a possible realization or scenario $\omega \in \Omega$ for the rest of the document, and therefore we write $t^\omega_j$ and $d^\omega_{jk}$ as realizations of these random parameters. 

The machines in set $\Set{M}$ are homogeneous, meaning they are all identical, and thus do not influence the different execution times $t_j^{\omega}$ and setup times $d_{jk}^{\omega}$ of each job $j, k \in \Set{J}$ for any realization $\omega \in \Omega$. Each machine has a time availability of $T$ and a maximum capacity of $B$ assigned jobs. The execution time and setup time of jobs assigned to a machine must meet the availability of time with a probability of at least $1 - \varepsilon$, with $\varepsilon \in (0, 1)$.

It is important to mention that this type of problem is not limited exclusively to solving CC-PMSP. The problem formulation is general enough to be used in other contexts with the same parameters and sets. For example, it can solve the Chance-Constrained Operating Room Scheduling (CC-ORS) problem where various surgical operations are assigned and sequenced in different operating rooms. Each of these operations has an execution time, which is not known with certainty because it may become complicated or may turn out easier than expected, and a stochastic setup time between operations, corresponding to the cleaning of the operating room and the preparations needed for the next intervention. The operating rooms also have a limited time and a maximum capacity of operations that can be assigned.

Another example that can be solved is the Chance-Constrained Vehicle Routing Problem (CC-VRP), where, instead of assigning and sequencing jobs to machines, different home deliveries are assigned and sequenced to different vehicles. Each of these deliveries has a stochastic execution time, which corresponds to the time it takes for the customer to receive the package, and a setup time that represents the travel time between customers. Each vehicle also has a limited time of use (for example, given the driver's work day) and a maximum capacity of orders to deliver.

What is interesting in these examples is the difference in the distribution of execution and setup times. For the CC-ORS problem, the execution time of a surgical operation is considerably longer than the setup time for the next operation. On the contrary, for the CC-VRP, the travel time between orders is significantly longer than the time it takes for the customer to receive the package. These two problems together allow us to perform more robust experiments, covering various configurations to evaluate the new formulations proposed in this research.

\subsection{Formulation with Chance Constraints}
\label{subsec:ChanceConstrFormulation}

One way to model the CC-PMSP is with the decision vector $\mathbf{x} \in \{0,1\}^{|\Set{J}|\times|\Set{M}|}$, where each variable $x_{jm}$ indicates whether job $j \in \Set{J}$ is assigned to machine $m \in \Set{M}$. The formulation is as follows:

\begin{subequations}
    \begin{alignat}{3}
        \max \qquad && \sum_{j\in \mathcal{J}} \sum_{m\in  \mathcal{M}} f_{j} x_{jm} \label{FO:2}\\ 
        \text{s.t.} \qquad && \sum_{m\in \mathcal{M}} x_{jm} &\leq 1 && \qquad \forall j \in \mathcal{J} \label{GenericConstr:1}\\
        && \sum_{j\in \mathcal{J}} x_{jm} &\leq B && \qquad \forall m \in \mathcal{M}\label{GenericConstr:2}\\
        && \mathbb{P}(\exists \mathbf{y} \in \Set{Y}(\mathbf{x}, \omega)) &\geq 1 - \varepsilon && \label{ChanceConstr:2}\\
        && x &\in \{0,1\}^{|\Set{J}|\times|\Set{M}|}
    \end{alignat}
\end{subequations}

The objective function \eqref{FO:2} maximizes the utilities of each job assigned to a machine. Constraints \eqref{GenericConstr:1} ensure that no job is assigned to multiple machines. Constraints \eqref{GenericConstr:2} impose that each machine is not assigned more jobs than its capacity limit.


Constraint \eqref{ChanceConstr:2} guarantees that $\Set{Y}(\mathbf{x}, \omega)$ is nonempty with probability at least $1 - \varepsilon$. In our problem, $\Set{Y}(\mathbf{x}, \omega)$ encodes sequencing the jobs assigned to each machine by $\mathbf{x}$ so that, given scenario $\omega$ of execution and setup times, the available time $T$ is not exceeded. The set $\Set{Y}(\mathbf{x}, \omega)$ is composed by vectors of variables $y_{jk}^{m\omega} \in \{0,1\}$ that indicate whether job $j \in \Set{J}\cup\{0\}$ precedes job $k \in \Set{J}\cup\{0\}$ on machine $m \in \Set{M}$ and a vector $\mathbf{u}^{m\omega}$ that ensures that the sequences of jobs assigned to all machines are valid. The set is defined as follows:
\begin{subequations}\label{SubproblemSet}
    \begin{align}
        \Set{Y}(\mathbf{x}, \omega) = \quad \{ \mathbf{y} \in &\R^{(|\Set{J}| +1) \times (|\Set{J}| +1) \times |\Set{M}|}: \nonumber\\
        &\sum_{k\in \mathcal{J} \cup \{0\}} y^{m\omega}_{jk} = x_{jm} & \forall j \in \mathcal{J}, m \in \mathcal{M} \label{PredecessorConstr} \\
        & \sum_{j\in \mathcal{J} \cup \{0\}} y^{m\omega}_{jk} = x_{km} &\forall k \in \mathcal{J}, m \in \mathcal{M} \label{SuccessorConstr}\\
        & \sum_{k\in \mathcal{J}\cup \{0\}} y^{m\omega}_{0k} = 1  &  \forall m \in \mathcal{M} \label{DummyConstr:1}\\
        & \sum_{j\in \mathcal{J}\cup \{0\}} y^{m\omega}_{j0} = 1 & \forall m \in \mathcal{M}\label{DummyConstr:2}\\
        & y_{jj}^{m\omega} = 0  & \forall j \in \mathcal{J}, m \in \mathcal{M} \label{LoopConstr}\\
        & \exists \mathbf{u}^{m\omega} \in \Set{S}_m(\mathbf{y}, \omega) & \forall m \in \mathcal{M} \label{SubToursConstr}\\
        & y_{jk}^{m\omega} \in \{0,1\} &\forall j, k \in \mathcal{J} \cup \{0\}, m \in \mathcal{M}
        \}
    \end{align}
\end{subequations}

The set employs the index $0$ as the dummy node used in scheduling problems to refer to the start and end of the sequence to be assigned to the jobs. Constraints \eqref{PredecessorConstr} and \eqref{SuccessorConstr} ensure that every job $j \in \Set{J}$ that has been assigned to machine $m \in \Set{M}$, is also assigned a position in the sequence of this machine. If $j$ has not been assigned to $m$, then the constraints ensure that this job is not in the machine's sequence. Constraints \eqref{DummyConstr:1} and \eqref{DummyConstr:2} ensure that the dummy node $0$ represents both the start and end of the job sequence, ensuring that there exists one and only one job associated with this start and end. Constraints \eqref{LoopConstr} limit that a job $j \in \Set{J}$ assigned to machine $m \in \Set{M}$ does not precede the same job $j$.

Constraints \eqref{SubToursConstr} guarantee a unique order $o$ in the sequence of jobs. We utilize the Single Commodity Flow (SCF) formulation \cite{scf} to model the scheduling problem. In our context, this formulation includes the following equations:
 \begin{subequations}\label{SCFEquations}
     \begin{align}
         &\Set{S}_m^{SCF}(\mathbf{y}, \omega) = \{\mathbf{u}^{m\omega} \in \R^{(|\Set{J}| +1) \times (|\Set{J}| +1)}:\nonumber\\
         &\sum_{k \in \mathcal{J} \cup \{0\}}u^{m\omega}_{jk} - \sum_{k \in \mathcal{J} \cup \{0\}} u^{m\omega}_{kj} = \sum_{k\in \mathcal{J} \cup \{0\}} y^{m\omega}_{jk} (t_j^{\omega} + d_{jk}^{\omega}) &\forall j \in \mathcal{J} \label{UDefinition}\\
         &u^{m\omega}_{jk} \leq T y^{m\omega}_{jk} & \forall j, k \in \mathcal{J}\cup \{0\}\label{LimitTimeConstr:2}\\
         &u^{m\omega}_{jk} \geq 0 & \forall j, k \in \mathcal{J} \cup \{0\}\}\label{NonNegativity:2}
     \end{align}
 \end{subequations}

The vector $\mathbf{u}^{m\omega} \in \mathbb{R}^{(|\Set{J}|+1)\times(|\Set{J}|+1)}$ acts as a time counter for the sequence of jobs assigned to each machine, so each variable $u_{jk}^{m\omega}$ indicates the time at which job $k \in \Set{J}\cup\{0\}$ can begin on machine $m \in \Set{M}$, taking into account the execution time of job $j \in \Set{J}\cup\{0\}$ and the setup time needed between them in scenario $\omega\in\Omega$.

Constraints \eqref{LimitTimeConstr:2} control that any pair of jobs $(j, k)$, assigned in that order to machine $m$, cannot exceed the time limit $T$ of said machine. Constraints \eqref{UDefinition} define the values of each variable within vector $\mathbf{u}^{m\omega}$, taking into account the accumulated time prior to job $j$, its execution time, and the setup time needed for the next job $k$ on machine $m$. Finally, \eqref{NonNegativity:2} imposes that the variables are non-negative, and together with \eqref{LimitTimeConstr:2}, force that the time of any pair of jobs $(j, k)$ that is not assigned to machine $m$ is equal to $0$.

We also explored the Miller-Tucker-Zemlin (MTZ) formulation \cite{mtz} as an alternative approach, but ultimately chose SCF due to its superior performance in preliminary experiments. The complete formulation and constraints of the MTZ formulation can be found in Appendix \ref{A:MTZ}.

\subsection{Formulation by Scenarios}
\label{subsec:FormulationScenarios}




From now on, we assume that $\Omega$ is finite.
Given that there are multiple scenarios $\omega \in \Omega$, we sum the probabilities associated with each scenario in which the scheduling problem constraints are satisfied and force this sum to be greater than or equal to the desired probability, $1 - \varepsilon$. For this, we introduce a new vector of variables $\mathbf{z} \in \{0,1\}^{|\Omega|}$ in which each $z^{\omega} \in \{0,1\}$ indicates whether the scheduling constraints are satisfied in scenario $\omega$. Thus, the chance constraint takes the following form:

\begin{subequations}
    $$\sum_{\omega \in \Omega} p^{\omega} z^{\omega} \geq 1 - \varepsilon$$
\end{subequations}
where $p^{\omega}$ represents the probability of occurrence of scenario $\omega \in \Omega$. We consider that all scenarios are equally probable, so each $p^{\omega}$ takes the value of $\frac{1}{|\Omega|}$.

The complete model with all deterministic constraints can be found in Appendix \ref{A:DeterministicModel}. Although a deterministic equivalent model provides an adequate representation of the problem, its complexity increases exponentially with the number of jobs, machines, and scenarios. Without proper decomposition, it becomes not possible to computationally solve this model. Therefore, in the following section, we present the approach we employed to avoid the limitations of the deterministic equivalent model.

\section{Decomposition Approach}
\label{sec:decomposition}

This section presents a decomposition for the CC-PMSP, in which the master problem assigns jobs to different machines, and the subproblems verify whether a job sequence exists that satisfies the time limit for each scenario and machine. Specifically, the master problem \eqref{MasterProblem} is defined as follows:
\begin{subequations}
    \begin{align}
        \max\ & \sum_{j\in \mathcal{J}} \sum_{m\in  \mathcal{M}} f_j x_{jm} \tag{\textsc{Master}} \label{MasterProblem}\\ 
        \text{s.t.} \ & \sum_{m\in \mathcal{M}} x_{jm} \leq 1 & \forall j \in \mathcal{J}\label{GenericConstr:3}\\
        & \sum_{j\in \mathcal{J}} x_{jm} \leq B &\forall m \in \mathcal{M}\label{GenericConstr:4}\\
        & \sum_{\omega \in \Omega} p^{\omega} z^{\omega} \geq 1 - \varepsilon \label{ChanceConstr:3}\\
        & x_{jm} \in \{0,1\} & \forall j \in \mathcal{J}, \ m \in \mathcal{M}\\
        & z^\omega \in \{0,1\} & \omega \in \Omega
    \end{align}
\end{subequations}

As detailed in Section \ref{sec:DescriptionFormulation}, variables $\mathbf{x}$ assign jobs to machines and $\mathbf{z}$ identify if a scenario satisfies the scheduling constraints. The objective function seeks to maximize the utilities of the assigned jobs. Constraints \eqref{GenericConstr:3} ensure that a job is not assigned to more than one machine, and constraints \eqref{GenericConstr:4} enforce that the machine limits are met. Lastly, \eqref{ChanceConstr:3} encodes the chance constraint as discussed in Section \ref{subsec:FormulationScenarios}.

As mentioned earlier, the subproblems are responsible for planning and ordering the jobs assigned to each machine. 
Note that the independence between machines allows for the decomposition to be not only in the set $\Omega$, but also in the set $\Set{M}$, so each subproblem $\mathcal{Q}(\hat{x}, \omega, m)$ depends on the candidate solution $\hat{x}$, the scenario $\omega \in \Omega$, and the machine $m \in \Set{M}$ to which the set of jobs was assigned. There are several possible ways to solve these subproblems, including the two DD-based formulations presented in Section \ref{sec:DDModels}. A possible IP formulation \eqref{SubProblemIP} for the subproblems is:
\begin{subequations}
    \begin{align}
        \mathcal{Q}(&\hat{x}, \omega, m) =  \max\  z \tag{SubIP} \label{SubProblemIP}\\
        \text{s.t.} \ & \sum_{k\in \mathcal{J} \cup \{0\}} y^{m\omega}_{jk} = \hat{x}_{jm} & \forall j \in \mathcal{J} \label{SubIP:Pertenencia1}\\
        &\sum_{j\in \mathcal{J} \cup \{0\}} y^{m\omega}_{jk} = \hat{x}_{km} & \forall k \in \mathcal{J}\label{SubIP:Pertenencia2}\\
        &\sum_{k\in \mathcal{J}} y^{m\omega}_{0k} = 1\\
        &\sum_{j\in \mathcal{J}} y^{m\omega}_{j0} = 1\\
        &y_{jj}^{m\omega} = 0 & \forall j \in \mathcal{J}\\
        &u_{jk}^{m\omega} - y^{m\omega}_{jk}(T - t_j^\omega) \leq M(1 - z) & \forall j, k \in \mathcal{J}\cup \{0\}\\
        &\sum_{k \in \mathcal{J} \cup \{0\}}u_{jk}^{m\omega} - \sum_{k \in \mathcal{J} \cup \{0\}} u_{kj}^{m\omega} = \sum_{k\in \mathcal{J} \cup \{0\}} y^{m\omega}_{jk} (t_j^\omega + d_{jk}^\omega) & \forall j \in \mathcal{J}\\
        &\mathbf{y}^{m \omega} \in \{0,1\}^{(|\mathcal{J}|+1)\times(|\mathcal{J}|+1)}\\
        &\mathbf{u}^{m\omega} \geq 0\\
        &z \in \{0,1\}
    \end{align}
\end{subequations}
The constraints are the same as in $\Set{Y}(\mathbf{x}, \omega)$ and $\Set{S}_m^{SCF}(\mathbf{y}, \omega)$. 

Given that we want to know if the subproblem is feasible, the objective function of $\mathcal{Q}(\hat{x}, \omega, m)$ indicates if it is feasible (i.e., $z = 1$) when there exists a job sequence that satisfies the limit $T$. If it is feasible, then the set of jobs that was assigned to machine $m \in \Set{M}$ is a possible valid configuration in scenario $\omega \in \Omega$. If it is not feasible, it means that there is no order for the assigned jobs that meets the machine's availability in scenario $\omega$, so one or more cuts should be added to the master problem.

The cuts we will use for the master problem have the following form:
\begin{equation}
    \sum_{j \in \bar{\Set{J}}^{h}} x_{jm} \leq |\bar{\Set{J}}^{h}| - z^\omega \qquad \forall h \in \Set{H}^\omega, \forall m \in \Set{M} \label{Cut}
\end{equation}
where $\Set{H}^\omega$ is the set of machines that, in scenario $\omega \in \Omega$, gave infeasible results with the candidate solution $\hat{x}$, and $\bar{\Set{J}}^{h}$ is a set of jobs that led to such infeasibility in machine $h \in \Set{M}$. Therefore, inequality  \eqref{Cut} ensures that the next candidate solution changes at least one job from the sets of jobs that were infeasible in past solutions. We use two types of cuts in our experimental evaluation: No-Good and IIS. In the case of a No-Good cut, the set $\bar{\Set{J}}^{h}$ contains all jobs assigned to $h$, that is, $\bar{\Set{J}}^{h} = \{j \in \Set{J}: \hat{x}_{jh} = 1\}$. On the other hand, the set $\bar{\Set{J}}^{h}$ for the IIS cuts covers the smallest combinations of jobs that make machine $h$ infeasible. Section \ref{sec:DDCuts} discusses the cuts in greater depth.

\subsection{Master Problem Valid Inequalities} \label{subsec:MasterProblemImprovements}
We now present a series of valid constraints that we add to the master problem \eqref{MasterProblem} to reduce the number of candidate solutions $\hat{x}$ that lead to equivalent and/or infeasible solutions.

\subsubsection{Symmetry Breaking Constraints} \label{subsubsec:SymmetryImprovements}
The CC-PMSP problem exhibits symmetries in its solutions due to the homogeneity of the machines. Thus, assigning the same set of jobs to one machine or another leads to equivalent solutions; however, the model does not recognize these equivalences, which can cause it to explore redundant solutions. In particular, the model can generate up to $|\Set{M}|!$ different solutions, where the sets of assigned jobs only exchange their machines. This redundancy introduces inefficiencies in the process of searching for the optimal solution, thereby increasing the resolution time without contributing to improvements in the final result.

To address this issue, we add the following symmetry-breaking constraints to \eqref{MasterProblem}:
\begin{subequations}\label{SymmetryConstr}
    \begin{alignat}{3}
        && x_{jm} &= 0 && \qquad \forall j \in \mathcal{J}, m \in \mathcal{M}: m > j \label{SymmetryConstr:1} \\
        && x_{j(m+1)} &\leq \sum_{k \in \mathcal{J}: k < j} x_{km} && \qquad \forall j \in \mathcal{J}, m \in \mathcal{M}: m < |\Set{M}| \label{SymmetryConstr:2}
     \end{alignat}
\end{subequations}

These constraints establish a lexicographic order in the sets of jobs assigned to the machines, identifying and excluding symmetric solutions. In particular, these constraints impose that the sets of jobs with lower indices are assigned to machines with smaller indices. Constraints \eqref{SymmetryConstr:1} guarantee that the index of a job is always less than the index of the machine to which it is assigned. Meanwhile, constraints \eqref{SymmetryConstr:2} maintain the lexicographic order by preventing a job with a lower index from being assigned to the previous machine, since if this happens, then the previous machine is lexicographically greater than the current machine and they should exchange jobs. We refer the reader to Appendix \ref{A:SymmetryExample} for more details on these symmetry-breaking constraints, where an example is provided that illustrates how these constraints effectively eliminate redundant solutions while preserving optimality.

\subsubsection{Improvements to Candidate Solutions}\label{subsubsec:CandidateImprovements}
A challenge identified during preliminary experiments was the large number of candidate solutions that do not satisfy the chance constraint. Each of these infeasible solutions forces us to solve the subproblems and add cuts, which implies a significant increase in the running time of our decomposition. To mitigate this effect, we introduce additional constraints to the master problem with more information about the subproblems. The objective of these constraints is to improve the quality of the generated candidate solutions, thereby reducing the occurrence of infeasible solutions and, consequently, enhancing the efficiency of the resolution process.

We implemented three types of constraints to augment the master problem, aiming to reduce the number of infeasible candidate solutions $\mathbf{\hat{x}}$.  The first type of constraints corresponds to a scenario relaxation approach, where we include one constraint per scenario and machine that represents a relaxation of the corresponding subproblem. 
These additional constraints ensure that the candidate solution $\mathbf{\hat{x}}$ satisfies the time limit of each machine for all scenarios $\{\omega \in \Omega: z^{\omega} = 1 \}$ with the most optimistic setup times, that is, 
\begin{equation}\label{eq:scenario-relaxation}
    \sum_{j \in \mathcal{J}} x_{jm} (t_j^{\omega} + \min_{k \in \mathcal{J}}\{ d_{jk}^{\omega} \}) \leq T + M(1 - z^{\omega}) \qquad \forall m \in \mathcal{M}, \omega \in \Omega 
\end{equation}

For each machine $m \in \Set{M}$ and for each scenario $\omega \in \Omega$, the constraint aggregates all execution times and the minimum setup times of the jobs assigned to machine $m$, ensuring that this sum does not exceed the machine's time limit. By considering the minimum of the setup times, this sum represents a lower bound on the time required to complete the entire sequence of jobs assigned to machine $m$ in scenario $\omega$. Consequently, this approach eliminates only those combinations that, even under the most favorable sequencing conditions, would be unable to satisfy the available time constraint.

We also consider two alternative constraint types: (i) include a small set of scenarios in the master problem, and (ii) introduce an optimistic scenario constraint. Unfortunately, preliminary results showed poor performance on these alternatives; thus, we briefly discussed them in Appendix \ref{A:ImprovementsToSolution}. We note that, although less effective for our specific problem instances, these alternative approaches may offer valuable insights for related CCP problems.

\section{DD-based Models for Sequencing Problem} \label{sec:DDModels}

This section introduces DD models utilized to solve the subproblem of our decomposition approach. We first formally define DDs for optimization problems and present the notation used throughout this paper. We then discuss how to use DDs to model scheduling problems and show the two proposed formulations for solving the subproblems of the CC-PMSP.

\subsection{Definition of DDs in Optimization} \label{subsec:DDOptimization}

To properly define DDs, consider the following optimization problem where $\textbf{x}$ is a $k$-dimensional vector of integer decision variables:
\begin{subequations}
    \begin{alignat}{2}
        \min \quad && f(\textbf{x}) \tag{P} \label{GeneralOptimizationProblem}\\
        \text{s.t.} \quad && \textbf{x} \in \Set{X} \nonumber
    \end{alignat}
\end{subequations}
Here $ f(\cdot)$ corresponds to a objective function and $\Set{X}$ represents the set of constraints that $\textbf{x}$ must satisfy to guarantee valid solutions. A DD for the problem \eqref{GeneralOptimizationProblem} is a directed acyclic graph $\Set{D} =  (\Set{N}, \Set{A})$, where its sets of nodes and arcs are represented by $\Set{N}$ and $\Set{A}$, respectively. The set of nodes $\Set{N}$ is partitioned into $k+1$ layers $\Set{L} = \{\Set{L}_1, ..., \Set{L}_{k+1}\}$, where $\bigcup_{i=1}^{k+1} \Set{L}_{i} = \Set{N}$ and $\Set{L}_i \cap \Set{L}_j = \emptyset$ for all $i \neq j$. The first and last layer always contain a single node called the root node ($\hat{r}$) and terminal node ($\hat{t}$), respectively.

The set of arcs is given by $\Set{A}$, where each arc $a=(n, n') \in \Set{A}$ connects an output node (i.e., parent node) $n \in \Set{N}$ with an input node (i.e., child node) $n' \in \Set{N}$. To refer to these nodes, we will use the notation $\textsf{out}(a) = n$ for the output node and $\textsf{in}(a) = n'$ for the input node. 
In general, arcs connect nodes belonging to consecutive layers of the diagram, that is, if $\textsf{out}(a) \in \Set{L}_i$, then $\textsf{in}(a) \in \Set{L}_{i+1}$. We will use $\delta^{in}(n)$ and $\delta^{out}(n)$ to represent the sets of incoming and outgoing arcs of a node $n \in \Set{N}$, respectively.

The DD of problem \eqref{GeneralOptimizationProblem} relates each of the layers $\Set{L}_i$, $i \in \{ 1,...,k \}$, with a decision variable $x_i$ and the outgoing arcs of that layer with the possible values of the variable. We will refer to $\textsf{val}(a)$ and $\textsf{c}(a)$ as the value assigned to arc $a \in \Set{A}$ and its cost for assigning that value, respectively. It is worth mentioning that no two arcs are leaving the same node that contain the same $\textsf{val}(a)$.

Each of the paths that can be formed from the root node to the terminal node is a solution to the problem \eqref{GeneralOptimizationProblem}. As we pass through the different layers through arcs $a \in \Set{A}$, values $\textsf{val}(a)$ are assigned to the various variables, and the costs $\textsf{c}(a)$ of the solution are added up. Therefore, to find the optimal solution, we need to find the path with the minimum cost. If the optimization problem were to maximize its objective function, we would need to find the path with the maximum cost.

To construct a DD for \eqref{GeneralOptimizationProblem}, we need a Dynamic Programming formulation of the problem. Specifically, we require four key components:  
\begin{enumerate*}[label=(\roman*)]
    \item a state space, where we denote the state of node $n \in \Set{N}$ as $\textsf{s}(n)$ and its domain of feasible arc values as $\textsf{dom}(n)$;
    \item a transition function $\textsf{tran}(\textsf{s}(n), \textsf{val}(a))$ that maps a state $\textsf{s}(n)$ and arc value $\textsf{val}(a)$ to a resulting state $\textsf{s}(n')$;
    \item a cost function that determines the cost $\textsf{c}(a)$ of an arc based on a state and arc value;
    \item and an initial state $\textsf{s}(\hat{r})$ from which the DD construction begins.
\end{enumerate*}

The DD construction is performed directly through the top-down algorithm construction approach developed by Bergman et al. \cite{bergman2014}. The specific algorithm adaptation for our problem is provided in Appendix \ref{A:TopDown}. We refer readers to the works of \cite{bergman2016,cire2013,vanhoeve2024} for a more comprehensive analysis of DD construction techniques.

\subsection{DD for Scheduling Problems} \label{subsec:DDScheduling}
We use two different DD-based models to solve the scheduling problem in our decomposition. In each iteration, each diagram must identify if there exists any order of the assigned jobs that meets the time availability of the machine and, if not, generate a cut for the master problem.

To model the order of jobs assigned to a machine, let us consider any machine $m \in \Set{M}$ and the set of jobs assigned to it in the master problem, denoted by $\Set{J}^m$. This set has a fixed number of jobs, and, thus, the DD employed to sequence the jobs and verify the feasibility of the subproblem has $|\Set{J}^m|+1$ layers. We use the decision vector $\mathbf{w} \in \Z_+^{|\Set{J}^m|}$ to assign each job to a position in the sequence, where the value of the variable $w_p$ indicates the index of the job assigned to position $p \in \{1,..,|\Set{J}^m|\}$.

Next, we present two DD-based formulations for solving the subproblems. The first formulation uses a DD with a linear cost function, in which the state of each node is composed of two components, one of which indicates the last job assigned. Due to this characteristic, we call this model DD-LJ, for the acronym DD-LastJob. The second formulation employs a DD that has a more complex cost function, designed to reduce the size of the diagram and simplify the state of the nodes. In this case, the state only contains the set of jobs already assigned. This model is called DD-JS, for the acronym DD-JobSet.

\subsubsection{Last Job DD-based Model} \label{subsubsec:DDLJ}
In DD-LJ, nodes use a state composed of two parts. The first part ($\textsf{s}_{set}$) consists of the jobs that have already been assigned to some position, while the second part ($\textsf{s}_{last}$) corresponds to the last job that was assigned a position. Formally, we have $\textsf{s}(n) := (\textsf{s}_{set}(n), \textsf{s}_{last}(n))$ with:
\begin{subequations}\label{DD-LJ:Estados}
    \begin{align}
        \textsf{s}_{set}(n) &:= \{w_{q}| q = 1, ..., (p-1)\}\\
        \textsf{s}_{last}(n) &:= w_{(p-1)}
    \end{align}
\end{subequations}
where $n \in \Set{L}_p$ is a node in $p$-th layer and $p \in \{1,..,|\Set{J}^m|\}$ is the position in the job order being assigned. The initial state is given by $\textsf{s}(\hat{r}) = (\{\}, -1)$, where the -1 is an indicator that there is no last job assigned. 

Given a node $n \in \Set{L}_p$ and an arc $a=(n,n')$, the state of node $n' \in \Set{L}_{p+1}$ consists of the same set of jobs that were already assigned in node $n$, along with the job represented by the arc. Then, the transition function for DD-LJ is:
\begin{equation} \label{DD-LJ:FuncionTransicion}
    \textsf{tran}(\textsf{s}(n), j) := (\textsf{s}_{set}(n') = \textsf{s}_{set}(n) \cup \{j\}, \textsf{s}_{last}(n') = j)
\end{equation}
where $j \in \Set{J}^m$ is a job such that $\textsf{val}(a) = j$.

For each node $n \in \Set{L}_p$, the domain of variable $w_p$ is given by $\textsf{dom}(n) = \Set{J}^m \setminus \textsf{s}_{set}(n)$. This means that the values that $w_p$ can take from node $n$ correspond to the jobs that are not present in $\textsf{s}_{set}(n)$.

The cost function for an arc $a \in \Set{A}$ and a scenario $\omega \in \Omega$ of the DD-LJ model is:
\begin{equation} \label{DD-LJ:FuncionCosto}
    \textsf{c}(a, \omega) = \begin{cases}
        t_{\textsf{val}(a)}^{\omega} & \text{if } \textsf{out}(a) = \hat{r} \\
        d_{\textsf{s}_{last}(\textsf{out}(a)),\textsf{val}(a)}^{\omega} + t_{\textsf{val}(a)}^{\omega} + d_{\textsf{val}(a), 0}^{\omega} & \text{if } \textsf{in}(a) = \hat{t}\\
        d_{\textsf{s}_{last}(\textsf{out}(a)),\textsf{val}(a)}^{\omega} + t_{\textsf{val}(a)}^{\omega} & \text{otherwise }
    \end{cases}\\
\end{equation}

That is, the cost of each arc $a \in \Set{A}$ is calculated as the setup time between the last job assigned in the ouput node, $\textsf{s}_{last}(\textsf{out}(a))$, and the job that the arc represents, $\textsf{val}(a)$, added to the execution time of the latter. The cost is different when they are arcs that enter the terminal node, as a setup time is added to represent leaving the machine ready for the next day. Additionally, the costs for arcs leaving the root node correspond only to the execution time of the first job since the machine is ready to be used. It is worth noting that this cost function is linear, as it only requires the information of the state of the ouput node $\textsf{out}(a)$ and the job arc $a$ represents to calculate $\textsf{c}(a)$ of any $a \in \Set{A}$. From now on, since the cost associated with an arc is a sum of times, we will call it the time of the arc.

\begin{figure}[ht]
    \begin{center}
        \includegraphics[scale=0.5]{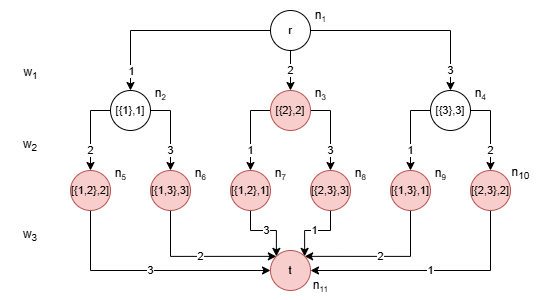} 
        \caption{DD-LJ structure for a subproblem with three jobs assigned to a machine}
        \label{Img:DD-LJ}
    \end{center}
\end{figure}

\begin{example}[DD-LJ]\label{Example:DD-LJ}
Let us take a machine $m \in \Set{M}$ to which the jobs $\Set{J}^m = \{1,2,3\}$ have been assigned. In Figure \ref{Img:DD-LJ}, we can see that the DD for this example has 3 layers, plus the final layer $\Set{L}_4$ that contains the terminal node. Each of the first three layers $\Set{L}_i$ is related to the variable $w_i$, defined at the beginning of the section. The numbers of each arc $a \in \Set{A}$ correspond to $\textsf{val}(a)$, which represent the job that was assigned in the current position.

In $\Set{L}_1$, we do not have any job assigned, so the domain of $w_1$ is $\textsf{dom}(\hat{r}) = \Set{J}^m$ and an arc is generated for each job. In $\Set{L}_2$, we already have the first job that will be processed by the machine assigned, so the domain of $w_2$ has two possible values and depends on the node $n \in \Set{L}_2$. If, for example, we arrive at $w_2$ passing through the node with state $(\{2\}, 2)$, it means that the job with index 2 has already been assigned in the first position of the sequence and that for the second position we can choose between the jobs with indices 1 and 3.

In $\Set{L}_3$, the last job that remains to be assigned is assigned, reaching the terminal node that already has all jobs assigned. For instance, if for $w_1$ and $w_2$ the jobs with indices 2 and 3 were chosen, respectively, then necessarily $w_3 = 1$. \hfill $\square$
\end{example}

\subsubsection{Job Set DD-based Model} \label{subsubsec:DDJS}
The second DD formulation used in the research was designed to create a smaller DD that represents the same set of solutions. For this, the states of the nodes of this DD are simply the jobs that have already been assigned a position, that is, $\textsf{s}(n) = \{w_{q}| q = 1, ..., (p-1)\}$, for all $n \in \Set{L}_p$. For the construction of this DD, we use an empty set as the initial state (i.e., $\textsf{s}(\hat{r}) = \{\}$). Then, the state of a node $n'$ with incoming arc $a=(n,n')$ corresponds to the set of jobs of the output node $n=\textsf{out}(a)$ and the job represented in such arc. So the transition function of this DD is
$$\textsf{tran}(\textsf{s}(n),j):= \textsf{s}(n) \cup \{j\}.$$

Similarly to DD-LJ, the values that $w_p$ can take from node $n \in \Set{L}_p$ correspond to the jobs that are not present in the state of the ouput node $\textsf{s}(n)$, i.e., $\textsf{dom}(n) = \Set{J}^m \setminus {\textsf{s}(n)}$.

Because the state of each node does not contain the index of the last job assigned, the cost function is more complex and dependeds on the information provided by arcs in the previous layer \cite{cire2013}.  Thus, the cost of each arc is calculated as follows:
\begin{equation} \label{FuncionCostoMDDCT}
    \textsf{c}(a, \omega) = \begin{cases}
        t_{\textsf{val}(a)}^{\omega} & \text{if } \textsf{out}(a) = \hat{r} \\
        \min\limits_{b \in \delta^{in}(n)} \{\textsf{c}(b,\omega) + d_{\textsf{val}(b),\textsf{val}(a)}^{\omega}\} + t_{\textsf{val}(a)}^{\omega} + d_{\textsf{val}(a),0}^{\omega} & \text{if } \textsf{in}(a) = \hat{t}\\
        \min\limits_{b \in \delta^{in}(n)} \{\textsf{c}(b,\omega) + d_{\textsf{val}(b),\textsf{val}(a)}^{\omega}\} + t_{\textsf{val}(a)}^{\omega} & \text{otherwise }
    \end{cases}\\
\end{equation}
where $a = (n,n')$. 

Note that to calculate the cost of each arc $a \in \Set{A}$, both the output node, $n = \textsf{out}(a)$, and the arcs entering it, $\delta^{in}(n)$, must be taken into account. Since the state does not have the last job assigned, it is necessary to recover the value of the setup time $d_{jk}^{\omega}$ of node $n$. For this, all the jobs that enter node $n$ are reviewed, and the one with the minimum time is selected, taking into consideration the job associated with arc $a$. This is where linearity is lost, as we must go back to the output node to calculate the cost of an arc. The difference when calculating the cost for arcs that meet $\textsf{out}(a) = \hat{r}$ or $\textsf{in}(a) = \hat{t}$ is the same as for the cost calculation in the DD-LJ model. Similarly to DD-LJ, since the cost associated with an arc is the sum of times, we call it the time of the arc.
\begin{figure}[ht]
    \begin{center}
        \includegraphics[scale=0.5]{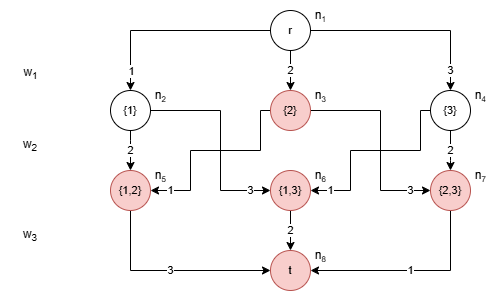} 
        \caption{DD-JS structure for a subproblem with three jobs assigned to a machine}
        \label{Img:DD-JS}
    \end{center}
\end{figure}
\begin{example}[DD-JS]
Let us consider the same case as Example \ref{Example:DD-LJ}, that is,  a machine with three jobs $\Set{J}^m= \{1,2,3\}$. The DD-JS can be seen in Figure \ref{Img:DD-JS}. We can notice that the first two layers, $\Set{L}_1$ and $\Set{L}_2$, are identical to those in Figure \ref{Img:DD-LJ}, but layer $\Set{L}_3$ changes. Since the state of a node does not depend on the last job assigned, in $\Set{L}_3$ we see that there are two ways to reach each node. For instance, if we take the node with state $\{1,2\}$, we can go through the node with state $\{1\}$ and through the node with state $\{2\}$. This not only simplifies the state, but also reduces the size of the DD. As the number of variables increases, this size reduction becomes greater. The red color of the nodes will be discussed in the next section. \hfill $\square$
\end{example}

\subsubsection{Network Flow DD-based Model} \label{subsubsec:NetworkFlow}
Lozano \& Smith \cite{lozano2018} propose a methodology that uses network flow models over a Binary Decision Diagram (BDD) to address a two-stage stochastic problem, in which the constraints of the second stage can be enabled or disabled depending on the candidate solution generated in the first stage. The decomposition used is a Benders' decomposition, in which the second stage is formulated to be solved with DD.

To adapt the network flow approach for our scheduling problem, we developed two types of DD: Capacitated Binary Decision Diagram (BDD-CAP) and Capacitated Multiple Decision Diagram (MDD-CAP). The former is the adaptation for scheduling problems of the BDD proposed by Lozano \& Smith \cite{lozano2018}, while the latter constitutes our attempt to enhance performance by reducing the size of the DD.

Unfortunately, neither of the network flow-based models demonstrated competitive computational performance for our purposes. Therefore, we have limited our discussion of these approaches in the main body of the paper, and present detailed notations, formulations, and explanations of each model in Appendix \ref{Asub:CapacitatedDD}. Additionally, Appendix \ref{Asub:NetworkFlowCuts} contains descriptions of the cuts implemented for our problem, adaptations of the strategies proposed by Lozano \& Smith \cite{lozano2018} to improve cut generation, and Appendix \ref{Asub:NetworkFlowResults} presents experimental results obtained through these approaches.

\section{DD-based Cuts} \label{sec:DDCuts}

We now detail how to use the DD models introduced in the previous section into our decomposition method and how to generate cuts from them. As previously mentioned, we utilize the DD models (i.e., DD-LJ or DD-JS) to solve each subproblem by verifying whether a job sequence exists with a time smaller than the time limit.

DDs have a great advantage: the structure of the DD only varies as a function of the number of layers that one wishes to construct. In the context of CC-PMSP, since the number of layers $|\Set{L}|$ depends on the number of jobs assigned to the machine whose feasibility is to be verified, it is possible to reuse the same DD structure for subproblems where the machine involved has the same number of assigned jobs. This property significantly reduces the time required to solve the subproblems. For a machine $m \in \Set{M}$, a DD with $|\Set{L}| = |\Set{J}^m| + 1$ layers is generated and then, for the remaining subproblems where the number of jobs assigned to a machine coincides with $|\Set{J}^m|$, the same DD can be used by only updating the cost function. As a result, instead of building or updating $|\Omega| \cdot |\Set{M}|$ subproblems, as happens in the \eqref{SubProblemIP} decomposition, we only need to construct at most $B$ diagrams.

Given a candidate solution $\hat{\mathbf{x}}$ from \eqref{MasterProblem}, we calculate the minimum sequencing time with one of our DDs for some specific machine and scenario. If this time does not exceed $T$, then the jobs assigned to the machine are feasible and a possible combination for the optimal solution. In case it is greater than $T$, then that job combination does not meet the time constraints for that scenario, and one or more cuts must be added to the master problem.

There are two types of cuts used in this work: No-Good and IIS. For DDs, the calculation of both cuts is relatively simple. In the case of No-Good cuts, the set $\bar{\Set{J}}^{h}$ is the set of jobs assigned to each machine $h \in \Set{M}$ that was infeasible in the candidate solution and the cut \eqref{Cut} is created. For example, for the DD in Figure \ref{Img:DD-JS}, if the arrangement were infeasible, the No-Good cut would simply take the complete set of jobs and generate the following cut for scenario $\omega \in \Omega$:
$$x_{1,m} + x_{2,m} + x_{3,m} \leq 3 - z^{\omega} \qquad \forall m \in \mathcal{M}$$

Additionally, DDs allow calculating several IIS per DD. As its name indicates, an IIS is a set of jobs that correspond to an infeasible assignment and for which any strict subset can lead to a feasible solution. Next, we detail the algorithms for finding IIS in DD-LJ and DD-JS.

\subsection{IIS Cuts with DD-LastJob Model} \label{subsec:DDLJIIS}

We use Algorithm \ref{alg:iisddlj} to calculate the IIS sets in the DD-LJ formulation. The shortest times to reach each node are stored in the dictionary $\texttt{possible\_IIS}$, where the keys represent sets of jobs, and the values correspond to the minimum scheduling times of those jobs. The algorithm considers that two different nodes $n, n' \in \Set{N}$ can share the same set of assigned jobs (i.e., $\textsf{s}_{set}(n) = \textsf{s}_{set}(n')$) and that there can be multiple paths to reach the same node. For example, in Figure \ref{Img:DD-LJ}, the nodes with states $\textsf{s}(n) = (\{1,3\},1)$ and $\textsf{s}(n') = (\{1,3\},3)$ share the same set $\textsf{s}_{set} = \{1,3\}$. To guarantee that $\{1,3\}$ is identified as an infeasible set of jobs, it is necessary to analyze both nodes.

\begin{algorithm}[htb]
    \caption{IIS Cut Creation for DD-LJ} \label{alg:iisddlj}
    \begin{algorithmic}[1]
        \Procedure{IISCut\_DD-LJ}{$\Set{D}$, T, $\omega$}
        \State $\texttt{time}(\hat{r})=0$
            \For{$l \in \{2,...,k+1\}$}
                \For{$n \in \Set{L}_l$}
                    \State $\texttt{time}(n) =$ Infinity
                    \If{$\textsf{s}_{set}(n) \notin$ \texttt{possible\_IIS}}
                        \State \texttt{possible\_IIS}[$\textsf{s}_{set}(n)$] = Infinity
                    \EndIf
                    \For{a $\in \delta^{in}(n)$}
                        \State $\texttt{current\_time} = \texttt{time}(\textsf{out}(a))$ $+$ $\textsf{c}(a,\omega)$
                        \If{\texttt{current\_time} $< \texttt{time}(n)$}
                            \State $\texttt{time}(n) = \texttt{current\_time}$
                            \If{$\texttt{time}(n) <$ \texttt{possible\_IIS}[$\textsf{s}_{set}(n)$]}
                                \State \texttt{possible\_IIS}[$\textsf{s}_{set}(n)$] = $\texttt{time}(n)$
                            \EndIf
                        \EndIf
                    \EndFor
                \EndFor
            \EndFor
            \State \texttt{IIS\_SET} = []
            \For{$\texttt{job\_set}$, $\texttt{min\_time} \in \texttt{possible\_IIS}$}
                \State \texttt{isIIS} = True
                \If{$\texttt{min\_time} > T$}
                    \For{\texttt{IIS} $\in$ \texttt{IIS\_SET}}
                        \If{\texttt{IIS} $ \subseteq \texttt{job\_set}$}
                            \State \texttt{isIIS} = False
                            \State \textbf{break}
                        \EndIf
                    \EndFor
                    \If{\texttt{isIIS} = True}
                        \State Append $\texttt{job\_set}$ to \texttt{IIS\_SET}
                    \EndIf
                \EndIf
            \EndFor
            \State \textbf{return} \texttt{IIS\_SET}
        \EndProcedure
    \end{algorithmic}
\end{algorithm}

Algorithm \ref{alg:iisddlj} iterates through all the nodes of the DD layer by layer, as shown in lines $\texttt{3}$ and $\texttt{4}$, while lines $\texttt{5}-\texttt{7}$ initialize the values of $\texttt{possible\_IIS}$ and $\texttt{time}(n)$. Then, in lines $\texttt{7-13}$, all the arcs entering the current node $\delta^{in}(n)$ are reviewed, the accumulated times are calculated using vector $\texttt{time}$ , and the minimum of them is stored in $\texttt{possible\_IIS}$.

Line $\texttt{14}$ defines the set that stores the cuts as $\texttt{IIS\_SET}$. In lines $\texttt{15}-\texttt{17}$, the algorithm reviews each of the tuples of \texttt{possible\_IIS}, which contains the job set and their corresponding minimum time, and compares them with the time availability $T$. If a time exceeds this limit, the algorithm evaluates whether the associated set of jobs is not only infeasible, but also irreducible (lines $\texttt{18-21}$). If that is the case, then the set of jobs is added to the \texttt{IIS} set (lines $\texttt{22-23}$). On the contrary, if it is not irreducible, it implies that at least one of its subsets already belongs to \texttt{IIS\_SET}. Note that all IISs of a particular size are included in \texttt{IIS\_SET} before iterating over the next size.


\begin{example}[IIS Algorithm for DD-LJ]
Let us consider Example \ref{Example:DD-LJ}, the DD in Figure \ref{Img:DD-LJ}, and the scenario $\omega \in \Omega$. Additionally, for this example, let us take the values $T = 5$, $t_1^{\omega} = 2$, $t_2^{\omega} = 6$, $t_3^{\omega} = 3$ and $d_{j,k}^{\omega} = 1$ for all $j\in \Set{J}^m$ and $k\in \Set{J}^m\cup\{0\}$.

The process of generating IIS cuts in DD-LJ begins with calculating the minimum time needed to reach each node. In the first layer, as we can see in \eqref{DD-LJ:FuncionCosto}, the times of the arcs $\delta^{out}(\hat{r})$ are $t^{\omega}_{\textsf{val}(a)}$, so their values are $2$, $6$ and $3$ for those reaching nodes $n_2$, $n_3$ and $n_4$, respectively. Since the root node has an accumulated time, \texttt{time}($\hat{r}$) in the algorithm, of $0$, then for this example we have that $\texttt{time}(n_2) = 2$, $\texttt{time}(n_3) = 6$ and $\texttt{time}(n_4) = 3$. The $\texttt{possible\_IIS}$ dictionary would add, in this layer, the following tuples: $(\{1\}, 2), (\{2\}, 6), (\{3\}, 3)$.

In \eqref{DD-LJ:FuncionCosto}, we can also see that the times of the arcs in the second layer are the sum of the setup time and the execution time. For the arc reaching node $n_5$, whose state is $(\{1,2\},2)$, its associated time is $d_{1,2}^{\omega} + t_2^{\omega} = 1 + 6 = 7$. Therefore, since \texttt{time}($n_2$) is the time required to reach the ouput node, $n_2$, plus the time of arc $(n_2, n_5)$, in this case, we have that $\texttt{time}(n_2) + \textsf{c}((n_2, n_5), \omega) = 2 + 7 = 9$. In the same way, the time required for each of the other nodes in the layer is calculated, which for this example are, $6$ for nodes $n_6$ and $n_9$, $9$ for node $n_7$ and $10$ for nodes $n_8$ and $n_{10}$. Since the $\texttt{possible\_IIS}$ dictionary stores the minimum time to reach each of the states, in this layer, the following tuples are added: $(\{1,2\}, 9), (\{1,3\}, 6), (\{2,3\}, 10)$.

Then, we have the terminal node in the last layer. Again, with \eqref{DD-LJ:FuncionCosto} we can calculate the time of each of the arcs reaching $\hat{t}$. The time of arc $(n_5, \hat{t})$ is calculated with $d_{2,3}^{\omega} + t_3^{\omega} + d_{3,0}^{\omega} = 1 + 3 + 1 = 5$ and if we add the time $\texttt{time}(n_5)$, we have that the accumulated time of $\hat{t}$ when arriving through arc $(n_5, \hat{t})$ is $5 + 9 = 14$. Similarly, when calculating the times of the other arcs, we have, for this case, that they are all $14$.

The accumulated time of the terminal node is $\texttt{time}(\hat{t}) = 14$ and since it is the only node with $\textsf{s}_{set} = \{1,2,3\}$, the $\texttt{possible\_IIS}$ dictionary adds the tuple: $(\{1,2,3\}, 14)$.

Finally, the algorithm evaluates the times of each set of jobs stored in $\texttt{possible\_IIS}$ to determine if any of them exceeds the time limit $T$, which in this case is $T = 5$. In Figure \ref{Img:DD-LJ}, the nodes marked in red represent the infeasible sets. It can be quickly verified that in the first layer, the only set that exceeds the availability of time is $\{2\}$. In $\Set{L}_2$, all sets exceed the value of $T$, but not all are IIS. The algorithm, specifically in line $\texttt{13}$, identifies that $\{2\}$ is a subset of the sets $\{1,2\}$ and $\{2,3\}$, which implies that these latter ones cannot be considered IIS. On the other hand, the set $\{1,3\}$ has a cost greater than $5$, and none of its subsets, $\{1\}$ and $\{3\}$, belongs to the IIS. Consequently, $\{1,3\}$ is added to the $\texttt{IIS\_SET}$. Finally, the set $\{1,2,3\}$ is analyzed, but since it contains subsets that are already IIS, it is not included in  $\texttt{IIS\_SET}$.

With this, the process of searching for IIS is concluded. As a result, the following cuts are added to the master problem:
\begin{subequations}
    \begin{alignat*}{3}
        && x_{2,m} &\leq 1 - z^{\omega} && \qquad \forall m \in \Set{M}\\
        && x_{1,m} + x_{3,m} &\leq 2 - z^{\omega} && \qquad \forall m \in \mathcal{M}
    \end{alignat*}
\end{subequations}
\hfill $\square$
\end{example}

\subsection{IIS Cuts with DD-JobSet Model} \label{subsec:DDJSIIS}
Algorithm \ref{alg:iisddjs} is used to calculate the IIS set of DD-JS cuts. The algorithm, in lines $\texttt{3}$ and $\texttt{4}$, traverses all the layers of the DD in order, evaluating each of its nodes. In lines $\texttt{5-9}$, the algorithm verifies if the state of the current node is a superset of an already found IIS. If it is, then the algorithm moves to the next node (lines \texttt{10} and \texttt{11}). If no IIS is a subset of the current node's state, then the algorithm in lines $\texttt{12-16}$ calculates the times of all incoming arcs of each node, storing the minimum time to arrive at the set of jobs. Finally, the minimum time is compared with the time limit $T$ (line $\texttt{17}$). This minimum value is used because, if it is less than or equal to $T$, there exists at least one sequence that satisfies the time availability constraint. In this case, the set of jobs is not infeasible and, therefore, cannot be part of the $\texttt{IIS\_SET}$. However, if the minimum time exceeds $T$, there is no feasible sequence for the set of jobs, thus fulfilling the definition of infeasibility. Therefore, the set is added to the $\texttt{IIS\_SET}$ in line $\texttt{18}$ of the algorithm.

\begin{algorithm}[tb]
    \caption{IIS Cut Creation for DD-JS} \label{alg:iisddjs}
    \begin{algorithmic}[1]
        \Procedure{IISCut\_DD-JS}{$\Set{D}$, T, $\omega$}
            \State \texttt{IIS\_SET} = []
            \For{$l \in \{2,...,k+1\}$}
                \For{$n \in \Set{L}_l$}
                    \State \texttt{superset} = False
                    \For{\texttt{IIS} $\in$ \texttt{IIS\_SET}}
                        \If{\texttt{IIS} $ \subseteq \textsf{s}(n)$}
                            \State \texttt{superset} = True
                            \State \textbf{break}
                        \EndIf
                    \EndFor
                    \If{\texttt{superset} = True}
                        \State \textbf{continue}
                    \EndIf
                    \State \texttt{time} = Infinity
                    \State \texttt{job\_set} = $\textsf{s}(n)$
                    \For{a $\in \delta^{in}(n)$}
                        \If {$\textsf{c}(a,\omega) <$ \texttt{time}}
                            \State \texttt{time} = $\textsf{c}(a,\omega)$
                        \EndIf
                    \EndFor
                    \If{\texttt{time} $>$ T}
                        \State Append \texttt{job\_set} to \texttt{IIS\_SET}
                    \EndIf
                \EndFor
            \EndFor
            \State \textbf{return} \texttt{IIS\_SET}
        \EndProcedure
    \end{algorithmic}
\end{algorithm}

As mentioned earlier, for a set of jobs to be considered an IIS, besides being infeasible it must be irreducible. That is, every strict subset of the set must have at least one feasible job sequence. This requirement is guaranteed in lines $\texttt{5-11}$ of the algorithm, where it is verified, for each node, if any set of the IIS is a subset of $\textsf{s}(n)$. In case it exists, the algorithm moves to the next node. In this way, we ensure that any set of jobs in $\texttt{IIS\_SET}$ does not have an infeasible subset.

Finally, once all the nodes of the DD have been evaluated, for each set of jobs identified as part of the IIS, the corresponding cuts are generated using equation \eqref{Cut}.

\begin{example}[IIS Algorithm for DD-JS]
Let us consider Example \ref{Example:DD-LJ}, the DD in Figure \ref{Img:DD-JS}, the scenario $\omega \in \Omega$, and let us assume that the times of the nodes painted red exceed $T$. Additionally, let us take the same values as in the previous example: $T = 5$, $t_1^{\omega} = 2$, $t_2^{\omega} = 6$, $t_3^{\omega} = 3$ and $d_{j,k}^{\omega} = 1$ for all $j\in \Set{J}^m$ and $k\in \Set{J}^m\cup\{0\}$.

The process of generating IIS cuts in DD-JS starts from the root node, which has time 0. Advancing sequentially, when evaluating the node with state $\{1\}$, whose time ($t_1^{\omega} = 2$) is less than $T = 5$, the progression continues. However, the node with state $\{2\}$ exceeds the value of $T$, since $t_2^{\omega} = 6 > T = 5$. Therefore, its state is incorporated into $\texttt{IIS\_SET}$. Continuing with the node with state $\{3\}$, the time $t_3^{\omega} = 3$ is also below $T$, so the algorithm advances to the next node.

The node with state $\{1,2\}$ contains an element previously included in $\texttt{IIS\_SET}$ (specifically $\{2\}$) and consequently, its revision is omitted. The same happens with the node $\{2,3\}$. Then, we have the node with state $\{1,3\}$. Following formula \eqref{FuncionCostoMDDCT}, we have that the cost of the arc is $\min\{(t_1^{\omega} + d_{1,3}^{\omega} + t_3^{\omega}), (t_3^{\omega} + d_{3,1}^{\omega} + t_1^{\omega})\} = 6$ which exceeds the value of $T=5$. Consequently, it is also added to $\texttt{IIS\_SET}$. And since the remaining nodes will always contain some element of $\texttt{IIS\_SET}$, the process of searching for IISs concludes. As a result, the following cuts are added to \ref{MasterProblem}:
\begin{subequations}
\begin{alignat*}{3}
    && x_{2,m} &\leq 1 - z^{\omega} && \qquad \forall m \in \Set{M}\\
    && x_{1,m} + x_{3,m} &\leq 2 - z^{\omega} && \qquad \forall m \in \mathcal{M}
\end{alignat*}
\end{subequations} \hfill $\square$
\end{example}

\subsection{Comparison of IIS Algorithms} \label{subsec:IISComparison}
Since the nodes of DD-LJ and DD-JS have different state structures, the algorithms for finding IIS differ quite a bit. The main difference is that for DD-LJ there can be several nodes that share the same $\textsf{s}_{set}(n)$. This characteristic increases the complexity of the algorithm, since, to determine the infeasibility of a set of jobs, it is necessary to examine all the nodes that contain it. For this reason, in Algorithm \ref{alg:iisddlj}, $\texttt{possible\_IIS}$ must be created, iterate through all nodes storing the minimum time to reach each set of jobs, and then compare these values with $T$ and verify that the state is irreducible.

Calculating IIS in Algorithm \ref{alg:iisddjs} is simpler because DD-JS nodes cannot share the same set of jobs, which corresponds to the complete state of the node, $\textsf{s}(n)$. This property allows us to directly determine whether the minimum time to reach the node is less than or greater than $T$ and, in case it is irreducible, identify it as an IIS. Also, an additional advantage arises in the implementation of Algorithm \ref{alg:iisddjs}: in certain cases, it is not necessary to traverse the entire DD-JS. If the algorithm finds a node whose state is an IIS, its descendant nodes can be automatically discarded, since they will be infeasible but not irreducible. This technique is not applicable in the implementation of Algorithm \ref{alg:iisddlj}, as the verification of irreducibility and infeasibility of the set of jobs is performed only after traversing all the nodes.

\section{Computational Experiments} \label{sec:Experiments}

We now present the experimental design used to evaluate the performance of the different models proposed throughout this research. We describe in detail the generated data sets, the parameters considered in the instances, and the computational environment used for the tests. Subsequently, we analyze the results obtained, examining both the overall performance of the models and their specific behavior under different configurations. This analysis enables us to understand the comparative advantages of DD-based approaches with traditional integer programming formulations.

\subsection{Problems and Instance Generation} \label{subsec:ProblemGeneration}
To rigorously evaluate the proposed models, we designed three data sets with different characteristics, inspired by real-world applications of the CC-PMSP. These sets represent scenarios where the distribution and proportion between execution times and setup times vary significantly, allowing us to examine the robustness of the formulations against different uncertainty structures. Below, we describe each of these sets in detail and the procedure used to ensure relevant mathematical properties such as the triangle inequality in setup times.

\paragraph{Operation Room Scheduling (ORS).} This data set is designed to address the CC-ORS problem presented in Section \ref{sec:DescriptionFormulation} and is characterized by operation execution times significantly greater than setup times between operations. To model this behavior, execution times follow a log-normal distribution with a mean of $\mu = 2$ and a standard deviation of $\sigma = 0.6$, while setup times are generated with parameters $\mu = 0.5$ and $\sigma = 0.15$ following the procedure described at the end of this section.
    
\paragraph{Vehicle Routing Problem (VRP).} Unlike the CC-ORS case, the CC-VRP problem, also presented in Section \ref{sec:DescriptionFormulation}, has a considerably smaller time to deliver the order to the customer than the travel time of a vehicle between destinations. To reflect this difference, we interchanged the means and deviations used in CC-ORS: execution times follow a log-normal distribution with a mean of $\mu = 0.5$ and a standard deviation of $\sigma = 0.15$, while setup times are generated with parameters $\mu = 2$ and $\sigma = 0.6$.
    
\paragraph{Equal.} Finally, we include a case where execution and setup times have the same mean and deviation. For this, we establish a mean of $\mu = 1.25$ and a standard deviation of $\sigma = 0.375$, maintaining the log-normal distribution for execution times and applying the procedure below to generate setup times.

To ensure that assigning more jobs to a machine does not decrease the total process time, we require that the setup times satisfy the triangle inequality. For this, we generate $|\Set{J}|+1$ points in a plane, maintaining an average distance of $\mu$ between them. Then, to generate uncertainty, in each scenario, a new point is generated from each original point by displacing it according to a random vector whose modulus follows a half-normal distribution with parameter $\sigma$ and whose angle is uniformly distributed in the interval $[0, 2\pi]$. The distances between each of these new points are the setup times for the scenario.

\subsection{Instance Parameters} \label{subsec:Parameters}

Our experiments consider nine configurations with different combinations of total number of jobs and machines, as well as different ratios $\frac{|\Set{J}|}{|\Set{M}|}$, to cover a wide range of situations and evaluate the robustness of the proposed models. These configurations allow us to analyze the behavior of each model under various conditions and understand its performance as a function of the relationship between jobs and machines. The chosen combinations are:
\begin{subequations}
\begin{align*}
    (|\Set{J}|,|\Set{M}|) \in \{&(60, 6), (80, 8), (100, 10), (72,6), (96,8),\\
    &(120,10), (84,6), (112,8), (140,10)\}
\end{align*}
\end{subequations}

To examine whether any model is favored by a larger number of jobs per machine, experiments with differentiated proportions are designed. The first three combinations present a ratio of 10 jobs per machine, the next three increase this ratio to 12 jobs per machine, and the last three reach 14 jobs per machine, thus allowing a more detailed evaluation of the impact of this relationship on model performance.

We consider the distribution of the data to define the value of the parameter $T$ (i.e., the time limit for each machine). In all data types, the sum of the averages of execution times and setup times is $2.5$. With this information, we approximated that on average each job will take $2.5$, counting its execution and the setup for the next job. Given that $B$ represents the maximum number of jobs that can be assigned to a machine, it was estimated that, on average, each machine would require a time of $T = B * 2.5$. We introduce the parameter $\textsc{dif}$ to regulate the difficulty of the problem, which allows us to adjust $T$ by either reducing it to increase complexity or increasing it to facilitate resolution. Thus, $T$ was parameterized as 
$T = 2.5B + \textsc{dif}*0.3$.

If the value of \textsc{dif} is too high, the time assigned to each machine increases, which trivializes the assignment and ordering of jobs. If, instead, \textsc{dif} is very negative, the time constraint becomes stricter, significantly increasing the difficulty of the problem and making it difficult to solve within the 20-minute limit. To ensure that the instances are neither too easy nor excessively complex, \textsc{dif} values are empirically selected to find intermediate solutions between these two extremes. The values that offered better results are: $\{0.2, 0.25, 0.3\}$ for ORS data, $\{-1, -0.95, -0.9\}$ for VRP data, and $\{-0.3, -0.25, -0.2\}$ for Equal data.

For the experiments, we calculate the big-M parameter $M$ by summing the execution time and the highest setup time for each job. Then, we selected the $B$ jobs with the highest values in this partial sum and added their resulting values to obtain the final estimate of $M$, that is:
$$M = \sum_{j \in \Set{J}_{topB}}\left(t_j + \max_{k \in \Set{J}: j \neq k}\{d_{ij}\}\right)$$
where $\Set{J}_{topB}$ are the $B$ jobs with the greatest duration.

For the remaining parameters, we consider $|\Omega| = 100$, $\varepsilon = 0.05$, and $B = \frac{|\Set{J}|}{|\Set{M}|}$ for all the experiments. These values were empirically determined with the objective of considering the largest possible instances that do not solve trivially and that, at the same time, present some probability of reaching an optimal solution within the 20-minute limit established for all experiments. A more detailed analysis of the choice of these parameters could be explored in future work. 

Overall, the experimental design encompasses three types of problems, nine combinations of $(|\Set{J}|,|\Set{M}|)$ per problem type, and three different values of $T$ for each one. Therefore, we have 27 distinct configurations for each problem type. We generate five instances for each configuration, resulting in a total of 135 instances per problem type and 405 instances in total.

\subsection{Computational Setup} \label{subsec:ComputationalSetup}

The algorithms are programmed in Python 3.11.6 using Gurobi 11.0.3 in the Visual Studio Code editor. All experiments were run on a 13th Gen Intel(R) Core(TM) i5-13600KF CPU with 16 GB of RAM. As mentioned earlier, each approach has a 20-minute time limit to find an optimal solution for an instance, which includes the time required to create all the models (i.e., IP and DDs models). The source code will be made available in a public repository upon publication.

To evaluate the new formulations presented in this work, we implement the DD-based decomposition using the DD-LJ and DD-JS models described in Section \ref{sec:DDModels}, using both No-Good and IIS cuts. In addition, we consider an IP-based decomposition where both \eqref{MasterProblem} and \eqref{SubProblemIP} are solved using Gurobi's technologies. This IP-based approach also considers a No-Good and an IIS alternative, where the IIS cuts are found using the tools provided by Gurobi. Additionally, preliminary experiments demonstrated that incorporating symmetric-breaking constraints \eqref{SymmetryConstr} and scenario relaxation constraints \eqref{eq:scenario-relaxation} was beneficial; therefore, these constraints are also included in the experiments below.

\subsection{Empirical Results for the DD-based Decomposition} \label{subsec:DDResults}

\begin{figure}
    \begin{center}
        \includegraphics[scale=0.45]{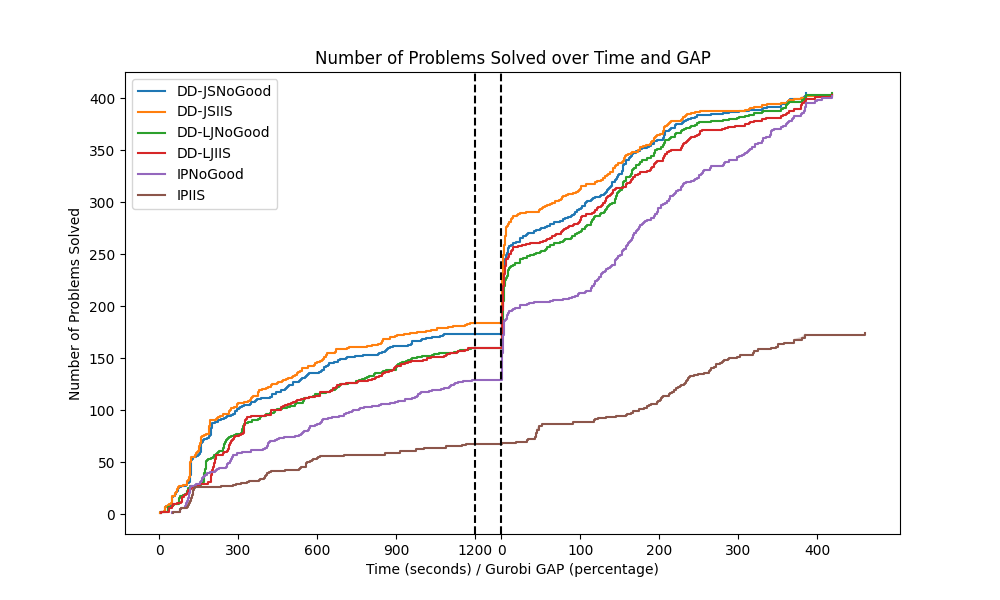} 
        \caption{Summary Profile Plot of All Experiments}
        \label{img:SummaryResults}
    \end{center}
\end{figure}

We now review the results of the 405 instances and the six different approaches: four DD-based decompositions (DD-LJ and DD-JS models using No-Good or IIS cuts) and two IP-based decompositions (i.e., using either No-Good or IIS cuts). Figure \ref{img:SummaryResults} presents a profile plot comparing the performance of the six approaches as a function of execution time (in seconds) and Gurobi's gap in percentage. On the X-axis, a combined scale is observed, where values to the left of the vertical line represent time in seconds, while values to the right correspond to Gurobi's gap in percentage. The Y-axis indicates the number of instances solved. Each colored line represents a different model  with the type of cut used.

The first black dotted vertical line marks the 20-minute time limit (i.e., 1200 seconds); thus, this first portion of the plot shows the number of instances solved to optimality. From here, we observe that the models that deliver the largest number of optimal solutions are the DD-JS models, then the DD-LJ models, and lastly the IP models. In particular, the DD-JS model with IIS cuts outperforms the other approach, solving 185 instances, while the best IP alternative (i.e., IP with No-Good cuts) only solves 129. Starting from the second black dotted line, the plot counts the instances that did not reach an optimal solution, starting from those that are closest to achieving the optimum (i.e., solutions with a small gap). We observe that the lines representing the DD-based decompositions consistently lie above those of the IP alternatives, indicating, in general, that the feasible solutions generated by the DD-based decompositions are closer to the optimum than those from IP models. 

\begin{table}
    \caption{Summary of results for all instances}
    \label{tab:SummaryResults}       
    \begin{tabular}{llrrr}
        \hline\noalign{\smallskip}
        Model & Cut & TotalTime & Gap & \#Optimal \\
        \noalign{\smallskip}\hline\noalign{\smallskip}
        DD-JS & IIS & 814.2 & inf & 184 \\
        DD-JS & No-Good & 835.3 & 0.6 & 173 \\
        DD-LJ & IIS & 896.5 & inf & 160 \\
        DD-LJ & No-Good & 902.5 & 0.7 & 160 \\
        IP & IIS & 1113.2 & inf & 64 \\
        IP & No-Good & 968.1 & inf & 129 \\
        \noalign{\smallskip}\hline
    \end{tabular}
\end{table}

Table \ref{tab:SummaryResults} supports the results presented in Figure \ref{img:SummaryResults}, where columns ``Model" and ``Cut" specify the model and the type of cut used, respectively. The ``TotalTime" column presents the average execution time, while ``Gap" shows the average value of the gap reported by Gurobi. Finally, the ``\#Optimal" column indicates the number of optimal solutions found within the 20-minute limit. We observe that decompositions based on DD-JS yield a greater number of optimal solutions and require less time compared to the DD-LJ alternative, which, in turn, outperforms the IP models in these metrics. Additionally, it is noteworthy that the DD alternatives with No-Good cuts are the only ones that always obtained a feasible solution, which explains why they are the only models without an infinite gap (i.e., ``inf'' symbol in the ``Gap'' column).

Table \ref{tab:SummaryResults2} presents filtered results that eliminate the instances where the models could not find a feasible solution. Here, the ``Model" and ``Cut" columns have the same meaning as in Table \ref{tab:SummaryResults}. The ``OptimalTimes" column represents the average execution time of the solutions that reached the optimum, while ``GapWithoutInf" shows us the average gap delivered by Gurobi for instances in which at least one feasible solution is found. Finally, the ``\#Inf" column indicates the number of instances where each approach could not find a feasible solution. Taking into account the information in Tables \ref{tab:SummaryResults} and \ref{tab:SummaryResults2}, in the case of DD-JS, we can see that, excluding the single instance where it failed to find a feasible solution, the model with IIS cut is, on average, superior to the No-Good, as it finds more optimal solutions in less time, and the suboptimal solutions have a slightly better gap. For the DD-LJ models, both are very similar, having the same number of optimal solutions and the same gap when excluding the instances that did not find a feasible solution. Finally, the No-Good cuts are far superior to IIS for the IP models.

\begin{table}[t]
    \caption{Summary of filtered results for all instances (i.e., only feasible solutions)}
    \label{tab:SummaryResults2}       
    \begin{tabular}{llrrr}
        \hline\noalign{\smallskip}
        Model & Cut & OptimalTimes & GapWithoutInf & \#Inf \\
        \noalign{\smallskip}\hline\noalign{\smallskip}
        DD-JS & IIS & 350.0 & 0.5 & 1 \\
        DD-JS & No-Good & 345.5 & 0.6 & 0 \\
        DD-LJ & IIS & 430.7 & 0.7 & 1 \\
        DD-LJ & No-Good & 420.1 & 0.7 & 0 \\
        IP & IIS & 384.0 & 1.8 & 229 \\
        IP & No-Good & 468.3 & 1.2 & 2 \\
        \noalign{\smallskip}\hline
    \end{tabular}
\end{table}

Table \ref{tab:SummaryResultsDetails} provides insights into why the DD-based approaches perform better compared to IP alternatives. As in the previous tables, the ``Model" and ``Cut" columns specify the model and the type of cut used, respectively. The ``\#Callback" and ``\#Cuts" columns indicate the average number of callbacks that each technique enters and the average cuts created, respectively. The ``ResolTime" column represents the sum of the times used in solving each subproblem, while ``ResolTimePerCB" shows the division between ``ResolTime" and ``\#Callback". The ``CreateCutTime" column indicates the sum of the time used in the creation of cuts, ``CreateSPTime" reflects the time dedicated to creating the subproblems, and finally, ``UpdateTime" indicates the time used to update the subproblems.

\begin{table}
    \caption{Detailed results of the decomposition approaches for all instances}
    \label{tab:SummaryResultsDetails}       
    \begin{tabular}{llrrrr}
        \hline\noalign{\smallskip}
        Model & Cut & \#Callback & \#Cuts & ResolTime & ResolTimePerCB \\
        \noalign{\smallskip}\hline\noalign{\smallskip}
        DD-JS & IIS & 95.8 & 20550.2 & 676.6 & 7.1 \\
        DD-JS & No-Good & 104.8 & 24134.6 & 692.8 & 6.6 \\
        DD-LJ & IIS & 76.5 & 16228.6 & 805.8 & 10.5 \\
        DD-LJ & No-Good & 92.4 & 21407.3 & 771.8 & 8.4 \\
        IP & IIS & 5.2 & 145.8 & 75.1 & 14.4 \\
        IP & No-Good & 43.0 & 5237.7 & 653.7 & 15.2 \\
        \noalign{\smallskip}\hline
    \end{tabular}
    \begin{tabular}{llrrrr}
        \hline\noalign{\smallskip}
        Model & Cut & CreateCutTime & CreateSPTime & UpdateTime \\
        \noalign{\smallskip}\hline\noalign{\smallskip}
        DD-JS & IIS & 0.5 & 0.3 & N/A \\
        DD-JS & No-Good & 0.7 & 0.3 & N/A  \\
        DD-LJ & IIS & 0.4 & 1.7 & N/A  \\
        DD-LJ & No-Good & 0.6 & 1.6 & N/A  \\
        IP & IIS & 746.5 & 246.9 & 21.2  \\
        IP & No-Good & 0.2 & 179.9 & 115.4  \\
        \noalign{\smallskip}\hline
    \end{tabular}
\end{table}

On average, DD techniques enter the callback significantly more times, allowing them to generate a greater number of cuts and, thus, find better solutions. Also, note that the poor performance of the IP with IIS cuts is primarily due to the time required to generate the cuts, which is four orders of magnitude larger than that of the No-Good alternative. Indeed, this explains why IP with IIS cuts enters fewer times to the callback and generates significantly fewer cuts than all other approaches. In contrast, the DD approaches take fractions of seconds to create the cuts, independently of whether it is a No-Good or an IIS cut. Furthermore, the DD models typically take less time to solve the subproblems than the IP models, with DD-JS being the fastest despite having a non-linear cost function. 




Lastly, we observe that as the number of machines increases, the problem becomes more complex. This is evident in Table \ref{tab:MachineResults}, where an increase in execution time and a decrease in the number of optimal solutions obtained are observed as the number of machines in the instance increases. Likewise, the gap tends to increase, except in cases where a feasible solution is not found.

\begin{table}[tb]
	\begin{center}
	\caption{Results by Amount of Machine}
    \begin{tabular}{llrrrr}
    \toprule
    Model & Cut & \#Machines & TotalTime & Gap & \#Optimal \\
    \midrule
    DD-JS & IIS & 6 & 644.1 & inf & 80 \\
    DD-JS & IIS & 8 & 812.3 & 0.3 & 65 \\
    DD-JS & IIS & 10 & 986.1 & 1.0 & 39 \\
    DD-JS & No-Good & 6 & 710.1 & 0.3 & 73 \\
    DD-JS & No-Good & 8 & 817.1 & 0.5 & 62 \\
    DD-JS & No-Good & 10 & 978.8 & 1.1 & 38 \\
    DD-LJ & IIS & 6 & 729.2 & inf & 76 \\
    DD-LJ & IIS & 8 & 903.7 & 0.6 & 51 \\
    DD-LJ & IIS & 10 & 1056.4 & 1.3 & 33 \\
    DD-LJ & No-Good & 6 & 758.8 & 0.4 & 67 \\
    DD-LJ & No-Good & 8 & 891.5 & 0.6 & 57 \\
    DD-LJ & No-Good & 10 & 1057.2 & 1.2 & 36 \\
    IP & IIS & 6 & 961.7 & inf & 36 \\
    IP & IIS & 8 & 1170.0 & inf & 15 \\
    IP & IIS & 10 & 1208.0 & inf & 13 \\
    IP & No-Good & 6 & 733.4 & 0.2 & 72 \\
    IP & No-Good & 8 & 1035.0 & 0.9 & 37 \\
    IP & No-Good & 10 & 1135.8 & inf & 20 \\
    \bottomrule
    \end{tabular}
    \label{tab:MachineResults}
	\end{center}
\end{table}

\subsection{Empirical Results for Each Problem Type} \label{subsec:DDResultsDataset}

To provide a more detailed analysis of model performance across different instance types, we present profile plots for three distinct datasets: VRP, ORS, and Equal. These figures allow us to visualize the performance differences between models for each specific dataset. Tables presenting the details of the results, broken down by problem type, are included in Appendix \ref{A:ResultsDataset}.

\begin{figure}[t]
    \begin{center}
        \includegraphics[scale=0.45]{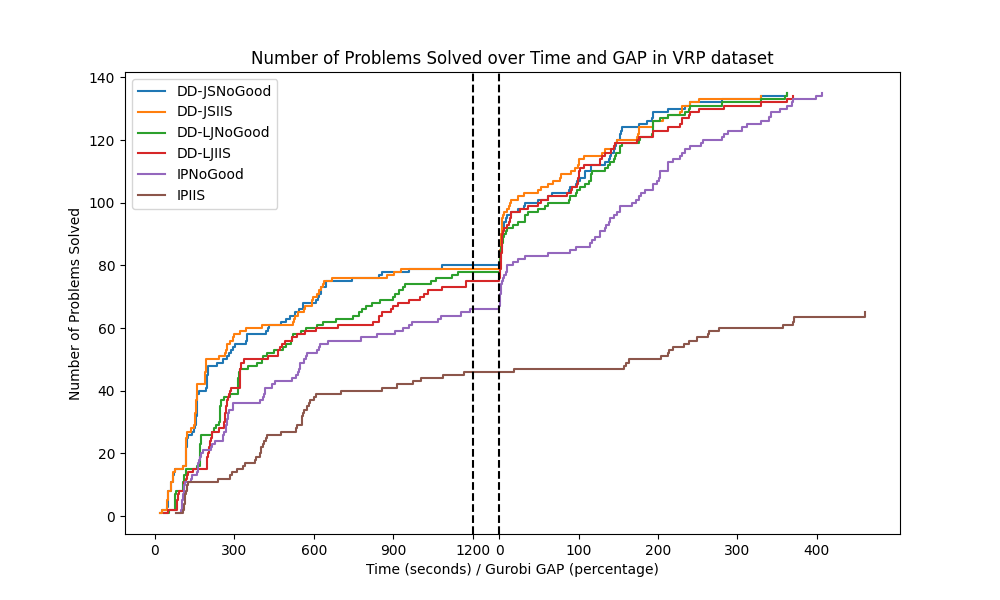} 
        \caption{Profile Plot for VRP data}
        \label{img:VRPProfilePlot}
    \end{center}
\end{figure}

Figure \ref{img:VRPProfilePlot}  displays the profile plot for the VRP dataset. In this dataset, DD-JS models with both cut types demonstrate superior performance, solving approximately 80 out of 135 instances to optimality before the time limit. The DD-LJ variants follow with 75 and 78 optimal solutions for IIS and No-Good, respectively. In contrast, IP models lag significantly, where the approach with No-Good cuts solves 65 instances, and the IIS cut alternative only 45. When considering feasible solutions, the DD-JS models continue to outperform the IP models, with their lines consistently positioned higher in the plot thus giving better sub-optimal solutions. As in Figure \ref{img:SummaryResults}, IP-IIS shows particularly poor performance in this dataset, with a flatter curve indicating both fewer optimal solutions and larger gaps in non-optimal solutions.

\begin{figure}[t]
    \begin{center}
        \includegraphics[scale=0.45]{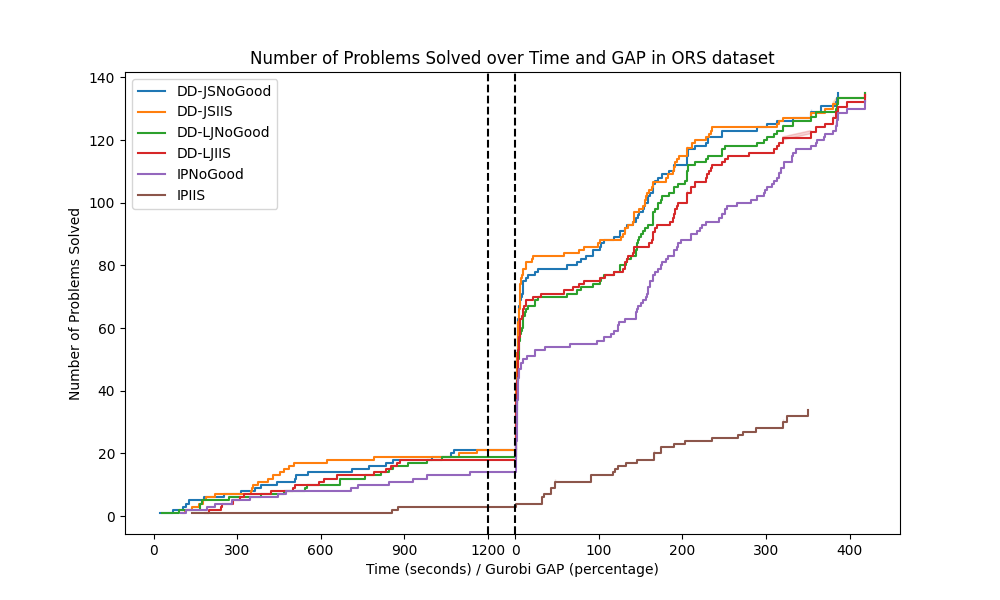} 
        \caption[ProfilePlotORS]{Profile Plot for ORS data}
        \label{img:ORSProfilePlot}
    \end{center}
\end{figure}

Figure \ref{img:ORSProfilePlot} presents the results for the ORS dataset. The performance pattern is similar to that shown in Figures \ref{img:SummaryResults} and \ref{img:VRPProfilePlot}, but the instances are more challenging to solve. DD-JS models achieve optimal solutions for 21 out of 135 instances, while DD-LJ and IP with No-Good cuts solve roughly 18 and 15 instances optimally, respectively. IP with IIS cuts performs particularly poorly in this dataset, finding only one optimal solution. In this plot, we can observe the most dramatic increase after the time limit, where we see that DD models consistently provide better solutions than IP models, with DD-JS notably returning better feasible solutions than DD-LJ.

\begin{figure}[t]
    \begin{center}
        \includegraphics[scale=0.45]{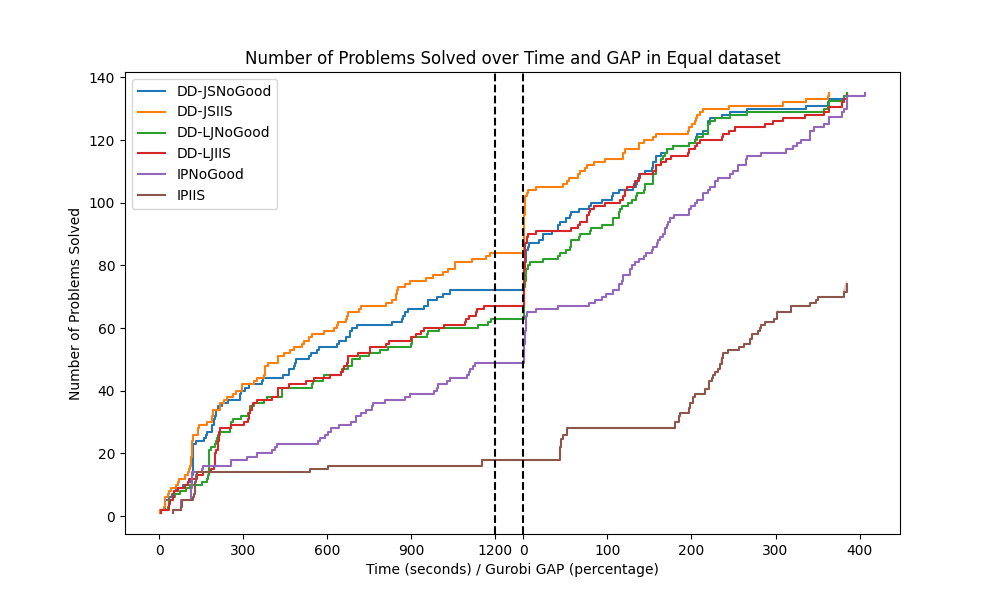} 
        \caption{Profile Plot for Equal data}
        \label{img:EqualProfilePlot}
    \end{center}
\end{figure}

Figure \ref{img:EqualProfilePlot} shows the performance for the Equal dataset, where the most notable feature is the exceptional performance of DD-JS with IIS cuts, which solves 84 out of 135 problems to optimality before the time limit. This is substantially better than other models, with DD-JS No-Good solving 72 instances, DD-LJ with IIS solving 67, DD-LJ with No-Good solving 63, and IP with No-Good solving only 49. IP with IIS shows significantly poorer performance in this dataset, solving fewer than 20 problems optimally. The difference between DD-JS with IIS and other models is particularly striking in this dataset, suggesting that this combination is especially well-suited for problems with equal distributions between realization and setup times. When analyzing the sub-optimal solutions, the DD-JS with IIS model maintains its superiority with a steeper curve, indicating consistently better gap values compared to all other approaches.

\subsection{Discussion}

The consistent performance patterns across these three datasets clearly demonstrate the superiority of DD-based models over IP models, with DD-JS variants generally outperforming DD-LJ variants across all problem types. While the general trend of DD-based superiority holds across all instance types, it is worth noting that the magnitude of this advantage varies, with particularly pronounced differences in the Equal dataset where DD-JS with IIS dramatically outperforms all other approaches. 

The results obtained demonstrate that DD-based decompositions significantly outperform IP alternatives in all evaluated datasets, both in the number of optimal solutions found and in the quality of the suboptimal solutions obtained. We identify two main reasons that explain the better performance of the DD-based techniques. First, the construction and update times of DD-based subproblems are considerably faster than those of subproblems formulated using IP. Second, the use of the shortest path algorithm proved to be quicker in solving the subproblems than the IP formulation. This translates to entering the callback more times, which allows for the generation of more cuts and, consequently, higher-quality solutions.

Likewise, the experiments demonstrate that the DD-JS model consistently outperforms the DD-LJ model. This is due to the same reasons mentioned above: the reduction in the size of the DD, in exchange for a more complex cost function, proved to be a beneficial strategy. A more compact DD allows for reducing construction times and shortest path resolution, which translates into being able to enter more callbacks, generate more cuts, and have better solutions.

In general terms, it cannot be categorically stated that IIS cuts outperform No-Good cuts. In the case of the Equal dataset, the approach with IIS cuts can, on average, find a greater number of optimal solutions in less time than models with No-Good cuts. However, in the other datacases, both types of cuts obtained practically identical results in most instances. We consider that this minimal differentiation could be attributed to the nature of the problems we generated, as in these experiments we observe that the IIS cuts are quite similar to No-Good cuts; that is, in most instances both cuts have the same or a similar infeasible set of jobs. This does not favor the IIS cuts because they are more computationally intensive, and our experiments do not clearly demonstrate whether their implementation is worthwhile in all problem types. Despite these findings, IIS cuts have similar or better performance than No-Good cuts in all three problem types and can be highly beneficial on problems where the underlying structure of the instances leads to more significant differences between the IIS and No-Good cuts.


\section{CONCLUSIONS} \label{sec:Conclusions}

This work explores different strategies based on Decision Diagrams (DD) to address scheduling problems with chance constraints, specifically the Chance-Constrained Parallel Machine Scheduling Problem (CC-PMSP). We formally define the problem and present a decomposition framework where the master problem assigns jobs to machines and the subproblems checks the feasibility of the joint chance constraints. Specifically, our decomposition involves a subproblem for each machine and scenario, where we need to verify whether the assigned jobs can be sequenced within the machine's time limit.

Given that the subproblems are NP-hard due to the presence of setup times, we explore the alternative of solving these subproblems using DDs. DDs are graphical structures that can represent the feasibility set of a problem. We leverage this property to build just a few decision diagrams (DDs) to describe the subproblems and reuse the DD graph to solve subproblems for different scenarios and machines (i.e., only changing the cost function). Moreover, we explore two DD models for these subproblems, where DD-LJ is larger but has a linear cost function, whereas DD-JS is smaller, but the cost computation is non-linear and thus more computationally expensive. 

We present novel algorithms to obtain No-Good and Irreducible Infeasible Set (IIS) cuts from the DDs, which can be potentially used for other applications. We also explore the use of Benders' cuts with DDs by adapting the techniques proposed by Lozano and Smith \cite{lozano2018}, but, unfortunately, with poor performance for our specific problem. Nonetheless, we believe that such methodology could be beneficial for other applications, and we encourage researchers to consider it.

We compare our DD-based alternatives with an Integer Programming (IP) based decomposition (i.e., using  IP technology to solve the subproblems and obtain No-Good or IIS cuts) across a wide range of instances coming from three different applications of the CC-PMPS (i.e., vehicle routing, operation room scheduling, and an ``equal time'' version of the PMSP). Our empirical results show a clear superiority of the DD-based alternatives, where our best alternative solves 55 more instances than the best IP decomposition (i.e., 184 vs. 124) and obtains significantly smaller gaps (i.e., 50\% vs. 120\%, on average). Our empirical evaluation presents the main findings in detail, including a discussion on why the DD-based alternatives outperform IP, a comparison of both DD models, and an analysis of using No-Good and IIS cuts.

In conclusion, this research presents a decomposition methodology based on DDs to solve scheduling problems with chance constraints, demonstrating high competitiveness. Although our work has focused specifically on the CC-PMSP, we believe that the methodology developed has significant potential to be applied to a broader range of problems. For example, the master problem could incorporate other types of constraints, focusing on entirely different domains. Likewise, the subproblems can be adapted to address other chance-constrained problems beyond scheduling, as long as they are amenable to DD technologies. This flexibility and adaptability constitute an additional contribution of our work, complementing the methodological advances already demonstrated and opening promising paths for future research at the intersection of DDs, chance constraints, and mathematical programming.


%
%

\bibliographystyle{spmpsci}      
\bibliography{Paper}   

\newpage
\appendix
\section{Miller-Tuckek-Zemlin Formulation} \label{A:MTZ}
The formulation of Miller-Tucker-Zemlin (MTZ) \cite{mtz}, in our context, includes the following constraints:

 \begin{subequations}\label{MTZEquations}
     \begin{align}
         \Set{S}_m^{MTZ}(\mathbf{y},\omega) &= \{\mathbf{v}^{m\omega}\in \R^{|\Set{J}| +1}: \nonumber  \\
         & v^{m\omega}_{j} - v^{m\omega}_{k} + t_j^{\omega} + d_{jk}^{\omega} \leq M (1 - y^{m\omega}_{jk}) &&\quad \forall j \in \mathcal{J}, k \in \Set{J}\cup \{0\} \label{VDefinition}\\
         & 0 \leq v^{m\omega}_{j} \leq T && \quad \forall j \in \mathcal{J}\cup \{0\} \label{NatureConstr:1} \}
     \end{align}
 \end{subequations}

The vector $\mathbf{v}^{m\omega} = (v^{m\omega}_0, ..., v^{m\omega}_{|\Set{J}|})$ functions as a time counter for the sequence of jobs assigned to machine $m$ in scenario $\omega$. Each variable $v^{m\omega}_{j}$ represents the moment at which job $j \in \Set{J}$ can begin on machine $m \in \Set{M}$. The variable associated with the dummy node, $v^{m\omega}_{0}$, indicates the total time needed to complete the entire sequence of jobs on that machine. The parameter $M$ is an upper bound on the total time required to complete any sequence of jobs.

Constraints \eqref{NatureConstr:1} limit the variables of vector $\mathbf{v}^{m\omega}$ within their possible values. In contrast, constraints \eqref{VDefinition} ensure the correct assignment of values to these variables, adding the corresponding execution and setup times. If a job $j \in \Set{J}$ precedes job $k \in \Set{J}$, then the right side of constraints \eqref{VDefinition} becomes zero and we force the time counter $v^{m\omega}_k$ to be greater than or equal to the counter $v^{m\omega}_j$ plus the corresponding execution and setup times. In case $j$ and $k$ are not together in the sequence, then we do not force any relationship between the time counters $v^{m\omega}_j$ and $v^{m\omega}_k$.

\section{Deterministic equivalent model} \label{A:DeterministicModel}

\begin{subequations}
    \begin{alignat}{3}
        \max \qquad && \sum_{j\in \mathcal{J}} \sum_{m\in  \mathcal{M}} u_{j} x_{jm} \tag{\textsc{DetEq}} \label{Deterministic}\\ 
        \text{s.t.} \qquad && \sum_{m\in \mathcal{M}} x_{jm} &\leq 1 && \qquad \forall j \in \mathcal{J} \\
        && \sum_{j\in \mathcal{J}} x_{jm} &\leq B && \qquad \forall m \in \mathcal{M}\\
        && \sum_{\omega \in \Omega} p^{\omega} z^{\omega} &\geq 1 - \varepsilon\\
        &&\sum_{k\in \mathcal{J} \cup \{0\}} y^{m\omega}_{jk} &= x_{jm} && \qquad \forall j \in \mathcal{J}, m \in \mathcal{M}, \omega \in \Omega\\
        && \sum_{j\in \mathcal{J} \cup \{0\}} y^{m\omega}_{jk} &= x_{km} && \qquad \forall k \in \mathcal{J}, m \in \mathcal{M}, \omega \in \Omega\\
        && \sum_{k\in \mathcal{J} \cup \{0\}} y^{m\omega}_{0k} &= 1 && \qquad \forall m \in \mathcal{M}, \omega \in \Omega \\
        && \sum_{j\in \mathcal{J} \cup \{0\}} y^{m\omega}_{j0} &= 1 && \qquad \forall m \in \mathcal{M}, \omega \in \Omega\\
        && y_{jj}^{m\omega} &= 0 && \qquad \forall j \in \mathcal{J}, m \in \mathcal{M}, \omega \in \Omega
    \end{alignat}
\end{subequations}

In case of opting to use the MTZ formulation, the following constraints are included:
\begin{subequations}
    \begin{alignat}{3}
        0 \leq v_{j}^{m\omega} &\leq T && \quad \forall j \in \mathcal{J}\cup \{0\}, m \in \mathcal{M}, \omega \in \Omega\label{MTZConstr:1}\\
        v_{j}^{m\omega} - v_{k}^{m\omega} + t_j^{\omega} + d_{jk}^{\omega} &\leq M (1 - y_{jk}^{m\omega}) && \quad \forall j \in \mathcal{J}, k \in \Set{J}\cup \{0\}, m \in \mathcal{M}, \omega \in \Omega \label{MTZConstr:2}
    \end{alignat}
\end{subequations}

In case of wanting to use the SCF formulation, the following constraints are applied:
\begin{subequations}
    \begin{alignat}{3}
        u_{jk}^{m\omega} &\leq T y_{jk}^{m\omega} &
        \begin{split}
            \qquad \forall j, k \in \mathcal{J}\cup \{0\},\\
            \qquad m \in \mathcal{M}, \omega \in \Omega
        \end{split}\label{SCFConstr:1}\\
        \sum_{k \in \mathcal{J} \cup \{0\}}u_{jk}^{m\omega} - \sum_{k \in \mathcal{J} \cup \{0\}} u_{kj}^{m\omega} &= \sum_{k\in \mathcal{J} \cup \{0\}} y_{jk}^{m\omega} (t_j^{\omega} + d_{jk}^{\omega}) &
        \begin{split}
            \qquad \forall j \in \mathcal{J}, m \in \mathcal{M}, \\
            \qquad \omega \in \Omega 
        \end{split}\label{SCFConstr:2}\\
        u_{jk}^{m\omega} &\geq 0 &
        \begin{split}
            \qquad \forall j, k \in \mathcal{J}, m \in \mathcal{M},\\
            \qquad \omega \in \Omega
        \end{split}\label{SCFConstr:3}
    \end{alignat}
\end{subequations}

\section{Symmetry Constraints Example} \label{A:SymmetryExample}
\begin{example}[Symmetry Constraints]
Let us take an example with 3 machines and 5 jobs, in which we have the following candidate solutions:
$$
\mathbf{x}^1 =
\begin{bmatrix}
    1 & 0 & 0\\
    0 & 1 & 0\\
    1 & 0 & 0\\
    0 & 0 & 1\\
    0 & 1 & 0\\
\end{bmatrix}
\qquad
\mathbf{x}^2 =
\begin{bmatrix}
    0 & 0 & 1\\
    0 & 1 & 0\\
    0 & 0 & 1\\
    1 & 0 & 0\\
    0 & 1 & 0\\
\end{bmatrix}
$$
We can verify that both solutions are symmetric, since in $\mathbf{x}^2$ the sets of jobs assigned to machines $1$ and $3$ of $\mathbf{x}^1$ were exchanged. But by adding constraints \eqref{SymmetryConstr}, $\mathbf{x}^2$ is no longer valid as it breaks both sets of constraints. Constraints \eqref{SymmetryConstr:1} are broken because $x_{1,3}^2 = 1$ and yet, the machine index is greater than the job index ($3 > 1$). But it also breaks constraints \eqref{SymmetryConstr:2}, because if we take $j = 2$ and $m = 1$ we have that $x_{2,2}^2 = 1$ and that $x_{1,1} = 0$, contradicting what the constraint asks us ($x_{2,2}^2 \leq x_{1,1}$). \hfill $\square$
\end{example}

\section{Other Improvements to Master Candidate Solutions}\label{A:ImprovementsToSolution}

\paragraph{Master Scenarios.} This strategy employs the first $\lfloor \varepsilon \cdot|\Omega| \rfloor + 1$ scenarios, which we will identify in the set $\Omega^{MS}$, and we add the following constraints to the master problem:
    \begin{subequations}
        \begin{alignat}{3}
            \sum_{\omega \in \Omega^{MS}} p^{\omega} z^{\omega} &\geq 1 - \varepsilon - \sum_{\omega \in \Omega \backslash \Omega^{MS}} p^{\omega} \label{MasterScenarios:1}\\
            \mathcal{Q}(x, \omega, m) &= 1 & \qquad \forall m \in \Set{M}, \omega \in \Omega^{MS} \label{MasterScenarios:2}
        \end{alignat}
    \end{subequations}
    In constraint \eqref{MasterScenarios:1} we seek that the candidate solution does not exceed the limit of unfulfilled scenarios in $\Omega^{MS}$. For example, consider $|\Omega| = 100$ and $\varepsilon = 0.1$, and therefore, $|\Omega^{MS}| = 100 \cdot 0.1 + 1 = 11$. The chance constraint \eqref{ChanceConstr:3} imposes that at most 10 scenarios do not comply with the scheduling restrictions \eqref{SubProblemIP}. Constraint \eqref{MasterScenarios:1} ensures that candidate solutions comply with \eqref{SubProblemIP} in at least one scenario $\omega \in \Omega^{MS}$, ensuring that this solution meets the requirements of \eqref{ChanceConstr:3} in at least the first scenarios. Constraints \eqref{MasterScenarios:2} are the same as in \eqref{SubProblemIP}, only that variable $x$ is not fixed.
    
\paragraph{Optimistic.} This strategy considers the most optimistic possible scenario $\omega^*$ and adds the following constraints to the master problem:
    \begin{subequations}
        \begin{alignat}{3}
            &\mathcal{Q}(x, \omega^*, m) = 1 && \qquad \forall m \in \Set{M} \label{Optimistic:1}
        \end{alignat}
    \end{subequations}
    The constraints are the same as in \eqref{SubProblemIP}, only that variable $x$ is not fixed and we force feasibility with $z = 1$. Additionally, it uses the optimistic scenario $\omega^*$ in which each job has an execution time of $t^{\omega^*}_j = \min\limits_{\omega \in \Omega}\{ t_{j}^{\omega} \}$ and setup times of $d^{\omega^*}_{jk} = \min\limits_{\omega \in \Omega}\{ d_{jk}^{\omega} \}$.

\section{Top-Down Algorithm} \label{A:TopDown}

\begin{algorithm}[htb]
\caption{DD Top-Down Construction} \label{alg:dd_topdown}
\begin{algorithmic}[1]
\Procedure{TopDownConstruction}{\textbf{x}}
\State Create DD $\Set{D}=(\Set{N}, \Set{A})$ with $|\textbf{x}| + 1 = k + 1$ empty layers.
\State Create the root node ($\hat{r}$) with the initial state.
\State Create the terminal node ($\hat{t}$) with a provisionally empty state.
\State The node $\hat{r}$ and the node $\hat{t}$ are added to layers $\Set{L}_1$ and $\Set{L}_{k + 1}$, respectively.
\For{$l \in \{1,...,k\}$}
\For{$n \in \Set{L}_l$}
\For{$v \in \textsf{dom}(n)$}
\State Create an arc $a$ with $\textsf{val}(a) = v$
\If{There exists a node $n' \in \Set{L}_{l+1}$ with $\textsf{s}(n') = \textsf{tran}(\textsf{s}(n), \textsf{val}(a))$}
\State Point arc $a$ to node $n'$, that is, $\textsf{in}(a) = n'$
\Else
\State Create node  $n'$ with $\textsf{s}(n') = \textsf{tran}(\textsf{s}(n), \textsf{val}(a))$
\State Add node $n'$ to layer $\Set{L}_{l+1}$
\State Point arc $a$ to $n'$, this is, $\textsf{in}(a) = n'$
\EndIf
\EndFor
\EndFor
\EndFor
\For{$n \in \Set{L}_{k}$}
\For{$v \in \textsf{dom}(n)$}
\State Create an arc $a$ with $\textsf{out}(a) = n$, $\textsf{val}(a) = v$ and $\textsf{in}(a) = \hat{t}$ 
\EndFor
\EndFor
\State \textbf{return} $\Set{D}$
\EndProcedure
\end{algorithmic}
\end{algorithm}

The DD construction is performed directly through Algorithm \ref{alg:dd_topdown}, which employs the Top-Down construction approach \cite{bergman2014}. The algorithm begins by initializing the DD structure, including the root and terminal nodes (lines \texttt{2-5}). Subsequently, the DD structure is traversed starting from the root node. For each node, an arc $a$ is generated corresponding to each of the values that its domain can assume (lines \texttt{7-9}).

In lines \texttt{10-15}, it verifies if there exists a node $n'$ with a state equivalent to the result of the transition function associated with arc $a$. If such a node exists, the arc is connected to node $n'$; otherwise, a new node is created and the arc is connected to this. Finally, in lines \texttt{16-17}, the connection between the nodes of the penultimate layer and the terminal node is established, completing the DD structure.

\section{Models Based on Network Flow} \label{A:NetworkFlowDD}




\subsection{Capacitated DD for Scheduling Problems} \label{Asub:CapacitatedDD}
To solve the network flow for our scheduling problem, we will present two types of DD: Capacitated Binary Decision Diagram (BDD-CAP) and Capacitated Multiple Decision Diagram (MDD-CAP). The first is the adaptation for scheduling problems of the BDD proposed by Lozano \& Smith \cite{lozano2018}, while the MDD is our contribution to reduce the size of the DD.

To adapt the capacitated arcs to our scheduling problem, two types of capacitated arcs were designed: assignment and non-assignment arcs. In the context of CC-PMSP, the requirements of these arcs relate to the assignment of jobs to a machine to allow flow in the diagram. In the case of assignment arcs, a specific job is required to be assigned to the machine for the flow to use that arc. On the other hand, in non-assignment arcs, one or more jobs are required not to be assigned to the machine to allow the flow. These requirements that each arc $a \in \Set{A}$ has to allow the flow will be identified with the set $\Set{U}_a$.

It is worth noting that, to use the capacitated arcs, the construction of both DDs used in this model includes all the jobs to be optimized, i.e., $\Set{J}$. This is because the capacitated arcs evaluate whether the assignment or non-assignment of a job to a machine could improve the result, so they need to evaluate all jobs. However, this could represent a disadvantage compared to the previously presented models, DD-LJ and DD-JS, which use only the jobs assigned to a certain machine, thus limiting the size of the layers to $|\Set{L}| \leq B$. In contrast, the layers of the DDs described in this section are bounded by $|\Set{J}|$, which implies a larger size.

\subsubsection{BDD-CAP Model} \label{Asubsub:BDDCAP}
This version of the model employs the BDD presented by Lozano \& Smith \cite{lozano2018}, whose structure differs from those used in DD-LJ and DD-JS. In this model, a decision vector $\mathbf{q} \in \{0,1\}^{|\Set{J}^m| \cdot |\Set{J}^m|}$ is defined, where each variable $q_{jp}$ indicates whether job $j$ has been assigned to position $p$ within the sequence. Consequently, for this type of DD, each layer$\Set{L}_{jp} \in \Set{L}$ is related to the variable $q_{jp}$, and for all nodes $n \in \Set{N}$, $\textsf{dom}(n) = \{0,1\}$. The value $1$ indicates that job $j$ is assigned to position $p$, while the value $0$ indicates that it will not be assigned to that position.

To construct the BDD-CAP, the same composite states defined in \eqref{DD-LJ:Estados} for DD-LJ are employed, as well as the same initial value $(\{\}, -1)$. However, the transition function differs, since in this case $\textsf{val}(a) \in \{0,1\}$ for all $a \in \Set{A}$. This implies that the values of the arcs no longer represent the job $j \in \Set{J}$ to which a position is to be assigned, but rather reflect a binary decision. The transition function of an arc $a = (n, n') \in \Set{A}$, where $n \in \Set{L}_{jp}$, is defined as:
\begin{equation} \label{BDD:FuncionTransicion}
    \textsf{tran}(\textsf{s}(n), \textsf{val}(a)) := \begin{cases}
        (\textsf{s}_{set}(n) \cup \{j\}, j) & \text{if } \textsf{val}(a) = 1 \\
        (\textsf{s}_{set}(n), \textsf{s}_{last}(n)) & \text{if } \textsf{val}(a) = 0
    \end{cases}\\
\end{equation}
It is important to note that, in case $\textsf{val}(a) = 1$, the job $j$ is used, which corresponds to the first index of the layer to which node $n$ belongs. On the other hand, when $\textsf{val}(a) = 0$, the state of the resulting node $n'$ is the same as the state of the parent node $n$, i.e., $\textsf{s}(n) = \textsf{s}(n')$.

We can differentiate three types of arcs in these DDs: non-capacitated arcs, assignment capacitated arcs, and non-assignment capacitated arcs. The sets of assignment capacitated arcs and non-assignment capacitated arcs will be identified with $\Set{A}_1$ and $\Set{A}_0$, respectively. The arcs that assign a job $j \in \Set{J}$ to position $p \in \{1,...,|\Set{J}|\}$, i.e., those arcs $a \in \Set{A}$ for which $\textsf{val}(a) = 1$ are assignment capacitated arcs and require that job $j$ be assigned to the machine. On the other hand, arcs that decide not to assign job $j$ to position $p$ and are directed to another node that is not the terminal, i.e., arcs $a \in \Set{A}$ with $\textsf{in}(a) \neq \hat{t}$ are non-capacitated arcs. And finally, arcs that do not assign job $j$ to position $p$ and are directed to the terminal node, ($\textsf{in}(a) = \hat{t}$) are non-assignment capacitated arcs and require that the jobs not assigned to a position are not assigned to the machine. These last arcs jump directly to the terminal node and are the only arcs that can skip layers.

The cost of each arc $a = (n, n')$ depends on whether it is an assignment arc, a non-assignment arc, or if it is not a capacitated arc, and is given by:
\begin{equation} \label{BDD:FuncionCosto}
    \textsf{c}(a,\omega) = \begin{cases}
        d_{\textsf{val}(a), 0}^{\omega } & \text{if } a \in \Set{A}_0 \\
        t_{\textsf{val}(a)}^{\omega} & \text{if } a \in \Set{A}_1 \text{ and } n' = \hat{r}\\
        d_{\textsf{s}_{last}(\textsf{out}(a)), \textsf{val}(a)}^{\omega} + t_{\textsf{val}(a)}^{\omega} + d_{\textsf{val}(a), 0}^{\omega} & \text{if } a \in \Set{A}_1 \text{ and } n = \hat{t}\\
        d_{\textsf{s}_{last}(\textsf{out}(a)), \textsf{val}(a)}^{\omega} + t_{\textsf{val}(a)}^{\omega} & \text{if } a \in \Set{A}_1\\
        0 & \text{if } a \notin \Set{A}_0 \cup \Set{A}_1
    \end{cases}\\
\end{equation}

We can see that the non-capacitated arcs have no cost. The non-assignment capacitated arcs, as they go directly to the terminal node, only add the setup time for node $\hat{t}$. Finally, the assignment capacitated arcs have the same logic explained in \eqref{DD-LJ:FuncionCosto} in which the execution time of the job and the corresponding setup times are added.

\begin{example}[BDD-CAP] \label{Example:BDD_CAP}
Let us consider that we want to assign three jobs $\Set{J} = \{1,2,3\}$ to two different machines. In Figure \ref{Img:BDD-CAP}, we can see the BDD-CAP that would be used to solve this problem. It is important to note that both subproblems associated with each machine would use the same BDD-CAP, which includes all jobs.

This time the DD has nine layers, since each of the three jobs can be in any of the three positions. Each of these layers $\Set{L}_{jp} \in \Set{L}$ is associated with the variable $q_{jp}$ and they decide whether job $j$ goes in position $p$ or not. The continuous arcs are those with $\textsf{val}(a) = 1$, while the discontinuous arcs are those with $\textsf{val}(a) = 0$. The black arcs are non-capacitated arcs, the green arcs are assignment capacitated arcs, and the red arcs are non-assignment capacitated arcs. The capacitated arcs have a set over the arrow that represents the set $\Set{U}_a$, that is, they represent their capacity. This set has the jobs that must be assigned or not assigned to the machine, depending on whether it is an assignment or non-assignment arc, for the flow to be able to pass through the arc.

\begin{figure}[ht]
    \begin{center}
        \includegraphics[scale=0.35]{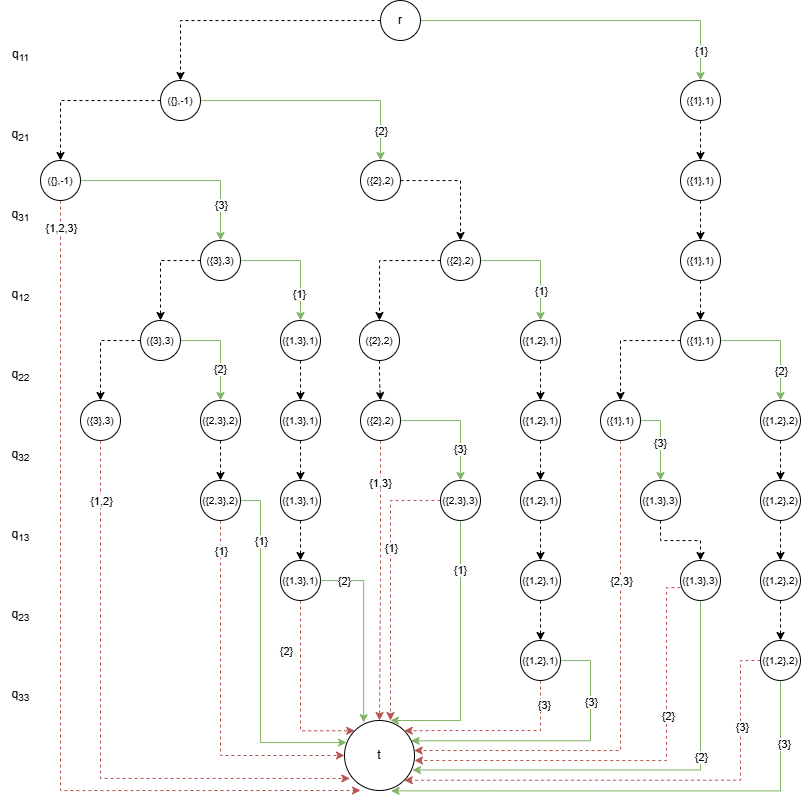} 
        \caption{BDD-CAP structure for a subproblem with three jobs assigned to a machine}
        \label{Img:BDD-CAP}
    \end{center}
\end{figure}

For example, let us consider the node that is in the center of the DD with state $(\{2\}, 2)$ that is just before assigning a value to the variable $q_{32}$, i.e., the DD must decide whether to include job $3$ in the second position of the sequence. From this node, two arcs emerge: one assignment capacitated (green) and another non-assignment capacitated (red). For the flow to traverse the green arc, since $\Set{U}_a = \{3\}$, job 3 must be assigned to the machine. Therefore, $\hat{\mathbf{x}}$ must be equal to $(1,1,1)$ or $(0,1,1)$ since job 2 has already been assigned. On the other hand, since the red arc is a non-assignment capacitated arc with $\Set{U}_a = \{1,3\}$, jobs 1 and 3 must not be assigned to the machine to allow the flow to pass through that arc, i.e., $\hat{\mathbf{x}} = (0,1,0)$. Finally, the non-capacitated arcs always allow the flow to pass. \hfill $\square$
\end{example}

\subsubsection{MDD-CAP Model} \label{Asubsub:MDDCAP}
The version of the model that uses an MDD-CAP is an adaptation of the BDD used by Lozano \& Smith \cite{lozano2018} and which we propose with the intention of reducing the size of the DD. The structure of the diagram is similar to the DDs proposed in Section \ref{sec:DDModels}, where the decision variables are the same as those used by DD-LJ and DD-JS. Specifically, we consider a vector $\mathbf{w} \in \Z_{+}^{|\Set{J}|}$ that assigns each job a position in the sequence, where the value of the variable $w_p$ indicates the index of the job assigned to position $p \in \{1,...,|\Set{J}|\}$ or takes the value $-1$ in case it is a non-assignment arc. Therefore, for each node $n \in \Set{N}$, the domain of the associated variable is $\textsf{dom}(n) = (\Set{J} \cup \{-1\})\setminus \textsf{s}_{set}(n)$.

It is worth noting that for this DD, the vector $\mathbf{w}$ is larger and the values of $p$ are bounded by the complete set of jobs and not by the set of jobs assigned to a machine $m \in \Set{M}$. This results in the MDD-CAP used in this model being considerably larger than the DD-LJ and DD-JS.

In this case, there are only two types of arcs and both are capacitated: the assignment arcs and the non-assignment arcs. These sets of arcs will maintain the same notation used in the BDD-CAP, i.e., $\Set{A}_1$ and $\Set{A}_0$, respectively. If the arc $a \in \Set{A}$ has $\textsf{val}(a) \neq -1$, then it is an assignment arc. On the contrary, if the arc $a$ has $\textsf{val}(a) = -1$, then it is a non-assignment arc.

To construct the MDD-CAP, the same states as in DD-LJ and BDD-CAP are used, defined in \eqref{DD-LJ:Estados}, and its initial value is also $(\{\}, -1)$. The transition function is similar to that presented in \eqref{BDD:FuncionTransicion}, only that the arcs can take more values.
\begin{equation} \label{MDD:FuncionTransicion}
    \textsf{tran}(\textsf{s}(n), \textsf{val}(a)) := \begin{cases}
        (\textsf{s}_{set}(n) \cup \{\textsf{val}(a)\}, \textsf{val}(a)) & \text{if } \textsf{val}(a) \neq -1 \\
        (\textsf{s}_{set}(n), \textsf{s}_{last}(n)) & \text{if } \textsf{val}(a) = -1
    \end{cases}\\
\end{equation}

The cost function is also similar to that presented in \eqref{BDD:FuncionCosto} for the BDD-CAP. We have that for all arcs $a = (n,n')$ and $\omega \in \Omega$:
\begin{equation} \label{MDD:FuncionCosto}
    \textsf{c}(a,\omega) = \begin{cases}
        d_{\textsf{val}(a), \hat{t}}^{\omega } & \text{if } a \in \Set{A}_0 \\
        t_{\textsf{val}(a)}^{\omega} & \text{if } a \in \Set{A}_1 \text{ and } n' = \hat{r}\\
        d_{\textsf{s}_{last}(\textsf{out}(a)), \textsf{val}(a)}^{\omega} + t_{\textsf{val}(a)}^{\omega} + d_{\textsf{val}(a), 0}^{\omega} & \text{if } a \in \Set{A}_1 \text{ and } n = \hat{t}\\
        d_{\textsf{s}_{last}(\textsf{out}(a)), \textsf{val}(a)}^{\omega} + t_{\textsf{val}(a)}^{\omega} & \text{if } a \in \Set{A}_1
    \end{cases}\\
\end{equation}

For the non-assignment capacitated arcs, the setup time for node $\hat{t}$ is added, since they always go to the terminal node. For the assignment capacitated arcs, the execution time, the corresponding setup time in case it is not the first assigned job, and the setup time towards $\hat{t}$ in case it is directed to the terminal node, must always be added.

\begin{figure}[ht]
    \begin{center}
        \includegraphics[scale=0.35]{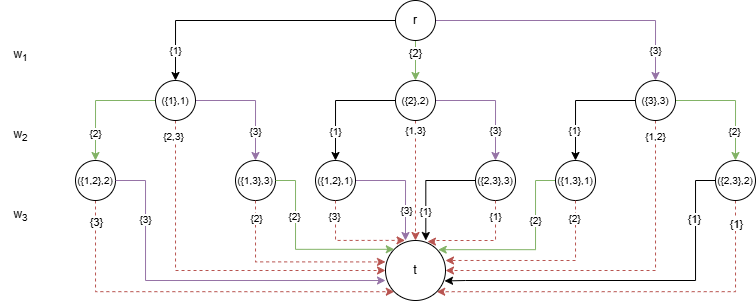} 
        \caption{MDD-CAP structure for a subproblem with three jobs assigned to a machine}
        \label{Img:MDD-CAP}
    \end{center}
\end{figure}

\begin{example}[MDD-CAP]
Let us consider the same set of jobs as in Example \ref{Example:BDD_CAP}. In Figure \ref{Img:MDD-CAP}, we can see the MDD-CAP that would be used to solve this problem. The MDD-CAP in the figure has three layers, one for each of the positions to which a job can be assigned, and each of these layers is associated with a variable $w_p$.

The continuous arcs are assignment arcs, while the discontinuous arcs are non-assignment arcs. Each of the colors of the arcs represents its value: black is $\textsf{val}(a) = 1$, green is $\textsf{val}(a) = 2$, purple is $\textsf{val}(a) = 3$, and red is $\textsf{val}(a) = -1$. The set over each of the arcs represents the set $\Set{U}_a$.

Let us consider the candidate solution $\hat{\mathbf{x}} = (1,0,1)$. From the root node, it is only possible to take the black and purple arcs, since the green arc is an assignment capacitated arc with capacity $\Set{U}_a = \{2\}$ and job 2 is not part of the solution.

From the node with state $(\{1\},1)$, the only available arc is the purple one. The green arc is still inaccessible for the same reason mentioned in the root node, while the red arc is a non-assignment capacitated arc with $\Set{U}_a = \{2,3\}$. However, since job 3 is assigned to the machine, the flow cannot transit through that arc.

Similarly, from the node with state $(\{3\},3)$, it is only possible to take the black arc. The non-assignment capacitated arc prevents the flow from passing, since its capacity includes job 1, which is assigned to the machine. In the same way, the assignment capacitated arc is not traversable due to its capacity including job 2.

Finally, in both paths, we arrive at a node with state $(\{1,3\},1)$ and another with state $(\{1,3\},3)$. In each case, there are two available arcs: one assignment capacitated and another non-assignment capacitated, both with $\Set{U}_a = \{2\}$. Since job 2 is not part of the solution, the only traversable arc is the red one.  \hfill $\square$
\end{example}
\subsection{Network Flow and Cuts} \label{Asub:NetworkFlowCuts}
Given one of the DDs described in the previous section, the network flow problem to solve our subproblem has the following form for a machine $m \in \Set{M}$ and scenario $\omega \in \Omega$:
\begin{subequations}
    \begin{alignat}{3}
        \text{min} && \sum_{a \in \Set{A}} \textsf{c}(a, \omega) \cdot v_{a} \label{FlujoEnRedes} \tag{\textsc{NetworkFlow}} \\
        \text{s.t.} &&  \sum_{a \in \delta^{out}(n)} v_{a} - \sum_{a \in \delta^{in}(n)} v_{a} &= \begin{cases}
            1\\
            0\\
            -1
        \end{cases} && \qquad \def\arraystretch{1.2}\begin{array}{@{}l}
            n = \hat{r} \\
            \forall n \in \Set{N} \backslash \{\hat{r}, \hat{t}\}\\
            n = \hat{t} 
        \end{array} \label{NetworkFlowConstr} \\
        && v_{a} &\leq \hat{x}_{qm} && \qquad \forall a \in \Set{A}_1, q \in \mathcal{U}_{a}\label{AsignationConstr} \\
        && v_{a} &\leq 1 - \hat{x}_{qm} && \qquad \forall a \in \Set{A}_0, q \in \mathcal{U}_{a} \label{NonAsignationConstr} \\
        && v_{a} &\geq 0 && \qquad \forall a \in \Set{A}
    \end{alignat}
\end{subequations}

The objective function seeks to minimize the times associated with the flow by passing from the root node $\hat{r}$ to the terminal node $\hat{t}$ in scenario $\omega$. The constraints defined in \eqref{NetworkFlowConstr} guarantee the logical consistency of the flow in the network and that the flow starts at $\hat{r}$ and ends at $\hat{t}$.

The set $\mathcal{U}_{a}$ contains the capacities of arc $a \in \Set{A}$, that is, it contains the jobs required for the flow to pass through the arc. In case $a \in \Set{A}_1$, this requires that the jobs $q \in \mathcal{U}_{a}$ be assigned to the machine. On the other hand, if $a \in \Set{A}_0$, the arc requires that the jobs $q \in \mathcal{U}_{a}$ not be assigned to the machine. Therefore, constraints \eqref{AsignationConstr} and \eqref{NonAsignationConstr} impose on the model to comply with the requirements of the assignment arcs and the non-assignment arcs. If a job $q \in \Set{U}_a$ of the capacity of arc $a \in \Set{A}_1$ is not assigned to machine $m \in \Set{M}$, i.e., $\hat{x}_{qm} = 0$, then in the corresponding constraints of \eqref{AsignationConstr} forces $v_a \leq 0$, preventing the flow from passing through that arc. Similarly, for an arc $a \in \Set{A}_0$, if a job $q \in \Set{U}_a$ is assigned to machine $m$, i.e., $\hat{x}_{qm} = 1$, then in the corresponding constraint of \eqref{NonAsignationConstr} forces $v_a \leq 1 - 1 = 0$, blocking the flow through the arc.

If we define $\pi$, $\alpha$ and $\beta$ as the dual variables of \eqref{NetworkFlowConstr}, \eqref{AsignationConstr} and \eqref{NonAsignationConstr}, respectively, the dual problem is of the following form:
\begin{subequations}
    \begin{align}
        \text{max} && \pi_{\hat{r}} - \pi_{\hat{t}} + \sum_{a \in A_1} \sum_{q \in \mathcal{U}_{a}} \hat{x}_{qm} &\cdot \alpha_{aq} + \sum_{a \in A_0} \sum_{q \in \mathcal{U}_{a}} (1 - \hat{x}_{qm}) \cdot \beta_{aq} \nonumber \tag{\textsc{DualNetworkFlow}} \label{FlujoEnRedesDual}\\
        \text{s.t.} && \pi_{n} - \pi_{n'} + \sum_{q \in \mathcal{U}_{a}} \alpha_{aq} &\leq \textsf{c}(a, \omega) && \forall a = (n, n') \in \Set{A}_1\\
        && \pi_{n} - \pi_{n'} + \sum_{q \in \mathcal{U}_{a}} \beta_{aq} &\leq \textsf{c}(a, \omega) && \forall a = (n, n') \in \Set{A}_0\\
        && \alpha &\leq 0\\
        && \beta &\leq 0
    \end{align}
\end{subequations}

Lozano \& Smith \cite{lozano2018} importantly note that, in general, there exist multiple valid dual solutions in \eqref{FlujoEnRedesDual}, which can affect the Benders cut in different ways. For this reason, they impose a specific interpretation of the dual variables that proved effective in their research. In their formulation, the optimal value $\hat{\pi}_{n}$ of the dual variable $\pi_{n}$ corresponds to the shortest path from node $n$ to the terminal node $\hat{t}$, which implies that $\hat{\pi}_{\hat{t}} = 0$. The optimal values $\hat{\alpha}_{aq}$ and $\hat{\beta}_{aq}$ of the dual variables $\alpha_{aq}$ and $\beta_{aq}$ represent reductions in the time of $\hat{\pi}_{\hat{r}}$ when the capacity $q \in \Set{U}_a$ of arc $a \in \Set{A}_1$ or $a \in \Set{A}_0$, respectively, is satisfied. Therefore, if the arc $a \in \Set{A}$ is included in $\hat{\pi}_{\hat{r}}$ or if its incorporation increases the time of the path, then it is fulfilled that $\hat{\alpha}_{aq} = 0$ or $\hat{\beta}_{aq} = 0$, depending on whether it is an assignment capacitated arc or a non-assignment arc.

To solve the master problem, we use optimality cuts of a Benders decomposition but adapted to take into account the chance constraint. Therefore, the cut we propose has the following form:
\begin{equation} \label{FlujoRedes:Corte}
    M (z^{\omega} - 1) + \hat{\pi}_{\hat{r}} + \sum_{a \in \Set{A}_1} \sum_{q \in \mathcal{U}_{a}} \hat{\alpha}_{aq} \cdot x_{qm} + \sum_{a \in \Set{A}_0} \sum_{q \in \mathcal{U}_{a}} \hat{\beta}_{aq} \cdot (1 - x_{qm}) \leq T
\end{equation}

As mentioned earlier, in the cut we have that $\hat{\pi}_{\hat{r}}$ represents the time of the shortest path in the current solution, while the terms $\hat{\alpha}_{aq}$ and $\hat{\beta}_{aq}$ correspond to possible reductions in that duration due to arc $a \in \Set{A}$, in case of adding job $q \in \Set{J}$ to machine $m \in \Set{M}$. If the sum of these terms, $\hat{\pi}_{\hat{r}} + \sum_{a \in \Set{A}_1} \sum_{q \in \mathcal{U}_{a}} \hat{\alpha}_{aq} \cdot x_{qm} + \sum_{a \in \Set{A}_0} \sum_{q \in \mathcal{U}_{a}} \hat{\beta}_{aq} \cdot (1 - x_{qm})$, is greater than $T$, then the time availability of the machine for that scenario is not met and therefore $z^{\omega} = 0$. This activates the parameter $M$ which is a sufficiently large number to omit the constraint. Otherwise, if the time availability is met, the use of parameter $M$ is not necessary since the inequality is fulfilled and consequently, $z^{\omega}$ can take the value of $1$.

\subsubsection{Improvements to the Cut} \label{Asubsub:NetworkFlowCutImprovements}
As mentioned at the beginning of the section, Lozano \& Smith \cite{lozano2018} propose two strategies to improve cut \eqref{FlujoRedes:Corte}. The first strategy takes advantage of the fact that in each layer only one arc can be taken. Therefore, adding all the decreases of all the arcs of a layer that benefit from adding a job is not a good bound. To adjust the bound, they define disjoint sets of arcs $\Gamma^{1}, ..., \Gamma^{D}$ such that $\Set{A}_1 = \bigcup_{d=1}^{D} \Gamma^{d}$ and $\Gamma^{k} \cap \Gamma^{l} = \emptyset$ for all $k,l \in \{1,...,D\}: k \neq l$. In these sets, every possible path from $\hat{r}$ to $\hat{t}$ uses at most one arc belonging to $\Gamma^d$ for all $d \in \{1,...,D\}$. The arcs that leave each layer in the DDs of our problem fulfill the definitions of $\Gamma$, so $D = |\Set{L}| - 1$. Applied to our case, the cut would be as follows:
\begin{equation} \label{FlujoRedes:CorteEstrategia1}
    M (z^{\omega} - 1) + \hat{\pi}_{\hat{r}} + \sum_{q = 1}^{|\mathcal{J}|} \sum_{d = 1}^D \hat{\gamma}_{qd} \cdot x_{qr} + \sum_{q = 1}^{|\mathcal{J}|} \hat{\delta}_{q} \cdot (1 - x_{qr}) \leq T
\end{equation}
where $\hat{\gamma}_{qd} = \min_{a \in \Gamma^{d}} \{ \hat{\alpha}_{aq} | q \in \mathcal{U}_{a} \}$ and $\hat{\delta}_q = \min  _{a \in \Set{A}_0} \{ \hat{\beta}_{aq} | q \in \mathcal{U}_{a} \}$. Therefore, this time we have the minimum time of the current solution ($\hat{\pi}_{\hat{r}}$) to which possibly at most the greatest reduction in time ($\hat{\gamma}_{qd}$) of each job $q \in \Set{J}$ in each layer $d \in \{1,...,D\}$ can be subtracted. Similarly, the same happens with $\hat{\delta}_q$ but with the non-assignment capacitated arcs.

This bound is tighter than the previous one because not all possible reductions of $\alpha$ and $\beta$ are taken into account, but only the minimum of them, that is, the one that can produce the greatest reduction in $\hat{\pi}_{\hat{r}}$. Note that $\delta$ does not need the sum of the disjoint sets of arcs because the non-assignment capacitated arcs all go to node $\hat{t}$ and, therefore, $\Set{A}_0$ is a set of disjoint arcs since no path from $\hat{r}$ to $\hat{t}$ can take two of these arcs.

The second strategy proposed by Lozano \& Smith \cite{lozano2018} consists of precomputing the optimal path for the network flow problem by forcing the flow to pass through each arc $\hat{a} \in \Set{A}$ such that $\mathcal{U}_{\hat{a}} \neq \emptyset$:
\begin{subequations}
\begin{align}
    l_{\hat{a}}^{\omega} = \quad &\text{min} \qquad & \sum_{a \in \Set{A}_1} \textsf{c}(a, \omega) \cdot v_{a} \nonumber\\
    & \text{s.t.} \qquad &  \sum_{i\in \mathcal{J}} x_{i} &\leq B\\
    && \sum_{a \in \delta^{out}(n)} v_{a} - \sum_{a \in \delta^{in}(n)} v_{a} &= \begin{cases}
            1\\
            0\\
            -1
        \end{cases} && \def\arraystretch{1.2}\begin{array}{@{}l}
            n = \hat{r} \\
            \forall n \in \Set{N} \backslash \{\hat{r}, \hat{t}\}\\
            n = \hat{t} 
        \end{array} \label{NetworkFlowConstr2} \\
    && v_{a} &\leq x_{q} && \forall a \in \Set{A}_1, q \in \mathcal{U}_{a}\\
    && v_{a} &\leq 1 - x_{q} && \forall a \in \Set{A}_0, q \in \mathcal{U}_{a}\\
    && v_{\hat{a}} &= 1\\
    && v_{a} &\geq 0 && \forall a \in \Set{A}\\
    && \mathbf{x} &\in \{0,1\}^{|\Set{J}|}
\end{align}
\end{subequations}

In this way, we obtain a new lower bound, since it represents the minimum possible time when using arc $a$. This allows us to adjust the bound in a better way. The resulting cut is as follows:
\begin{equation} \label{FlujoRedes:CorteEstrategia2}
    M (z^{\omega} - 1) + \hat{\pi}_{\hat{r}} + \sum_{q = 1}^{|\mathcal{J}|} \sum_{d = 1}^D \textrm{\underline{$\hat{\gamma}$}}_{qd}^{\omega} \cdot x_{qr} + \sum_{q = 1}^{|\mathcal{J}|} \textrm{\underline{$\hat{\delta}$}}_{q}^{\omega} \cdot (1 - x_{qr}) \leq T
\end{equation}
where
\[
\textrm{\underline{$\hat{\gamma}$}}_{qd}^{\omega} = \min_{a \in \Gamma^{d}} \{ \textrm{\underline{$\hat{\alpha}$}}_{aq}^{\omega} \mid q \in \mathcal{U}_{a} \}, \quad 
\textrm{\underline{$\hat{\alpha}$}}_{aq}^{\omega} = \max\{\hat{\alpha}_{aq}, (l^{\omega}_{a} - \hat{\pi}_{\hat{r}})^-\},
\]
\[
\textrm{\underline{$\hat{\delta}$}}_{q}^{\omega} = \min_{a \in \Set{A}_0} \{ \textrm{\underline{$\hat{\beta}$}}_{aq}^{\omega} \mid q \in \mathcal{U}_{a} \}, \quad 
\textrm{\underline{$\hat{\beta}$}}_{aq}^{\omega} = \max\{\hat{\beta}_{aq}, (l^{\omega}_{a} - \hat{\pi}_{\hat{r}})^-\}.
\]

The intuition behind this cut is to take advantage of the best lower bounds obtained in the previous strategies and incorporate them into the formulation. It is important to note that the coefficients $\textrm{\underline{$\hat{\gamma}$}}$ and $\textrm{\underline{$\hat{\delta}$}}$ are obtained through a minimization of the maximum between the lower bounds, thus guaranteeing the validity of the inequality. In particular, the terms $\textrm{\underline{$\hat{\alpha}$}}$ and $\textrm{\underline{$\hat{\beta}$}}$ adjust the lower bound, and the minimization ensures that the constraint remains valid in all cases.

\subsection{Results of the Network Flow Model} \label{Asub:NetworkFlowResults}

The reason why the results of the Network Flow Model are not included in the previous results is that it was not possible to create decision diagrams with the dimensions required by these models. Unlike the decision diagrams of the MDD Model, which have a depth equal to the number of jobs assigned to the machine, the decision diagrams of the Network Flow Model have a depth equal to the total number of jobs of the complete problem. Therefore, the depth can be up to $|\mathcal{M}|$ times greater.

If we calculate the nodes of the decision diagrams of the DD-LJ model, which will always have equal or greater quantity than the nodes of the DD-JS, we find that the maximum number of nodes is $B!$. However, when calculating the number of nodes of the Network Flow Models, we obtain that they are, at most, $|\mathcal{J}|!$ for the MDD-CAP and $2^{|\mathcal{J}|\cdot|\mathcal{J}|}$ for the BDD-CAP. Therefore, the difference in sizes between the decision diagrams is considerable.

For these models, we had to lower the number of scenarios to $|\Omega| = 50$ and lower the number of jobs and machines $(|\Set{J}|,|\Set{M}|) \in \{(6, 2), (8, 2)\}$. The rest of the parameters remain exactly the same.
\begin{table}
    \caption{Summary of Unfiltered Network Flow Results}
    \label{tab:SummaryResultsLozano}       
    \begin{tabular}{llrrr}
        \hline\noalign{\smallskip}
        Model & Strategy & TotalTime & Gap & \#Optimal \\
        \noalign{\smallskip}\hline\noalign{\smallskip}
        BDD-CAP & - & 645.1 & inf & 45 \\
        BDD-CAP & 1 & 629.0 & inf & 45 \\
        MDD-CAP & - & 595.6 & inf & 71 \\
        MDD-CAP & 1 & 605.6 & inf & 68 \\
        \noalign{\smallskip}\hline
    \end{tabular}
\end{table}

\begin{table}
    \caption{Summary of Filtered Network Flow Results}
    \label{tab:SummaryResultsLozano2}       
    \begin{tabular}{llrrr}
        \hline\noalign{\smallskip}
        Model & Strategy & OptimalTimes & GapWithoutInf & \#Inf \\
        \noalign{\smallskip}\hline\noalign{\smallskip}
        BDD-CAP & - & 79.0 & - & 45 \\
        BDD-CAP & 1 & 70.6 & - & 44 \\
        MDD-CAP & - & 439.1 & 0.1 & 5 \\
        MDD-CAP & 1 & 403.6 & 0.1 & 6 \\
        \noalign{\smallskip}\hline
    \end{tabular}
\end{table}
Tables \ref{tab:SummaryResultsLozano} and \ref{tab:SummaryResultsLozano2} have the same columns as Tables \ref{tab:SummaryResults} and \ref{tab:SummaryResults2}. The only difference is that, as the Network Flow models do not have No-Good and IIS cuts, that column was replaced by the strategies proposed by Lozano \& Smith \cite{lozano2018}. It is important to mention that under the second strategy in these models, no results could be seen even by lowering the parameters to the minimum because the optimization of $l_{\hat{a}}^{\omega}$ for all arcs took more than the 20-minute time limit.

As can be observed in the tables, despite the considerable reduction in the number of jobs, machines, and scenarios, these models show limited performance. The MDD-CAP model manages to find 71 optimal solutions without a specific strategy and 68 with strategy 1, while the BDD-CAP reaches only 45 optimal solutions in both cases, out of a total of 90 instances evaluated for each configuration.

Table \ref{tab:SummaryResultsLozano2} reveals that the average time to reach optimal solutions is significantly lower for BDD-CAP (79.0 and 70.6 seconds) compared to MDD-CAP (439.1 and 403.6 seconds). However, the BDD-CAP model failed to find feasible solutions in 45 and 44 instances, respectively, which represents practically 50\% of the cases. In contrast, MDD-CAP only failed to find feasible solutions in 5 and 6 instances, respectively, also achieving an average gap of 0.1\% for sub-optimal solutions.

When contrasting these results with those obtained by the DD-JS and DD-LJ models presented earlier, the difference is notable. While the network flow models required drastically reducing the size of the instances to just 6-8 jobs and 2 machines with 50 scenarios, the DD-JS and DD-LJ models were able to address problems with up to 140 jobs, 10 machines, and 100 scenarios. According to Table \ref{tab:SummaryResults}, DD-JS managed to find 184 and 173 optimal solutions (with IIS and No-Good cuts respectively), while DD-LJ reached 160 optimal solutions with both types of cuts, out of a total of 243 instances per configuration.

\section{Experiment Results by data sets} \label{A:ResultsDataset}
\begin{table}[H]
	\begin{center}
	\caption{Results by data sets}
    \begin{tabular}{lllrrr}
    \toprule
    Model & Cut & Dataset & \#Problems & \#Optimal & TotalTime \\
    \midrule
    DD-JS & IIS & Equal & 135 & 84 & 709.9 \\
    DD-JS & IIS & ORS & 135 & 21 & 1080.3 \\
    DD-JS & IIS & VRP & 135 & 79 & 652.4 \\
    DD-JS & NoGood & Equal & 135 & 72 & 758.8 \\
    DD-JS & NoGood & ORS & 135 & 21 & 1090.1 \\
    DD-JS & NoGood & VRP & 135 & 80 & 657.1 \\
    DD-LJ & IIS & Equal & 135 & 67 & 825.2 \\
    DD-LJ & IIS & ORS & 135 & 18 & 1109.6 \\
    DD-LJ & IIS & VRP & 135 & 75 & 754.6 \\
    DD-LJ & NoGood & Equal & 135 & 63 & 844.6 \\
    DD-LJ & NoGood & ORS & 135 & 19 & 1122.7 \\
    DD-LJ & NoGood & VRP & 135 & 78 & 740.3 \\
    IP & IIS & Equal & 135 & 17 & 1121.1 \\
    IP & IIS & ORS & 135 & 1 & 1237.5 \\
    IP & IIS & VRP & 135 & 46 & 981.1 \\
    IP & NoGood & Equal & 135 & 49 & 957.6 \\
    IP & NoGood & ORS & 135 & 14 & 1132.7 \\
    IP & NoGood & VRP & 135 & 66 & 813.9 \\
    \bottomrule
    \end{tabular}
    \label{tab:ResultsDataset}
	\end{center}
\end{table}

\subsection{VRP} \label{Asub:VRPResults}

\begin{table}[H]
	\begin{center}
	\caption{Summary of VRP Unfiltered Results}
    \begin{tabular}{llrrrrr}
    \toprule
    Model & Cut & \#Problems & \#Optimal & TotalTime & Gap \\
    \midrule
    DD-JS & IIS & 135 & 79 & 652.3986 & inf \\
    DD-JS & No-Good & 135 & 80 & 657.1437 & 0.4160 \\
    DD-LJ & IIS & 135 & 75 & 754.5817 & inf \\
    DD-LJ & No-Good & 135 & 78 & 740.2876 & 0.4813 \\
    IP & IIS & 135 & 46 & 981.0909 & inf \\
    IP & No-Good & 135 & 66 & 813.8522 & 0.8578 \\
    \bottomrule
    \end{tabular}
	\label{tab:VRPSummaryValues}
	\end{center}
\end{table}

\begin{table}[H]
	\begin{center}
	\caption{Summary of VRP Filtered Results}
    \begin{tabular}{llrrrr}
    \toprule
    Model & Cut & OptimalTimes & GapWithoutInf & \#Inf \\
    \midrule
    DD-JS & IIS & 263.8 & 0.4 & 1 \\
    DD-JS & No-Good & 283.5 & 0.4 & 0 \\
    DD-LJ & IIS & 397.8 & 0.4 & 1 \\
    DD-LJ & No-Good & 393.0 & 0.5 & 0 \\
    IP & IIS & 448.9 & 0.7 & 70 \\
    IP & No-Good & 408.5 & 0.9 & 0 \\
    \bottomrule
    \end{tabular}
    \label{tab:VRPSummaryValues2}
	\end{center}
\end{table}

\begin{table}[ht]
	\begin{center}
	\caption{Detailed VRP Results}
    \begin{tabular}{llrrrr}
    \toprule
    Model & Cut & \#Callback & \#Cuts & ResolTime & ResolTimePerCB \\
    \midrule
    DD-JS & IIS & 97.9 & 28134.4 & 500.1 & 5.1 \\
    DD-JS & No-Good & 99.3 & 32223.9 & 504.1 & 5.1 \\
    DD-LJ & IIS & 78.0 & 21939.0 & 656.0 & 8.4 \\
    DD-LJ & No-Good & 88.4 & 29043.5 & 601.0 & 6.8 \\
    IP & IIS & 3.5 & 99.0 & 77.6 & 22.3 \\
    IP & No-Good & 33.7 & 5847.4 & 527.3 & 15.7 \\
    \bottomrule
    \end{tabular}
    \begin{tabular}{llrrrr}
    \toprule
    Model & Cut & CreateCutTime & CreateSPTime & UpdateTime & ResolTime \\
    \midrule
    DD-JS & IIS & 0.7 & 0.3 & N/A & 500.1 \\
    DD-JS & No-Good & 0.9 & 0.2 & N/A & 504.1 \\
    DD-LJ & IIS & 0.6 & 1.4 & N/A & 656.0 \\
    DD-LJ & No-Good & 0.8 & 1.4 & N/A & 601.0 \\
    IP & IIS & 613.5 & 241.4 & 23.4 & 77.6 \\
    IP & No-Good & 0.2 & 179.2 & 92.7 & 527.3 \\
    \bottomrule
    \end{tabular}
    \label{tab:VRPSummaryResultsDetails}
	\end{center}
\end{table}

\subsection{ORS} \label{Asub:ORSResults}

\begin{table}[H]
	\begin{center}
	\caption{Summary of ORS Unfiltered Results}
    \begin{tabular}{llrrrrr}
    \toprule
    Model & Cut & \#Problems & \#Optimal & TotalTime & Gap \\
    \midrule
    DD-JS & IIS & 135 & 21 & 1080.3 & 0.8 \\
    DD-JS & No-Good & 135 & 21 & 1090.1 & 0.8 \\
    DD-LJ & IIS & 135 & 18 & 1109.6 & 1.1 \\
    DD-LJ & No-Good & 135 & 19 & 1122.7 & 1.0 \\
    IP & IIS & 135 & 1 & 1237.5 & inf \\
    IP & No-Good & 135 & 14 & 1132.7 & inf \\
    \bottomrule
    \end{tabular}
	\label{tab:ORSSummaryValues}
	\end{center}
\end{table}

\begin{table}[H]
	\begin{center}
	\caption{Summary of ORS Filtered Results}
    \begin{tabular}{llrrrr}
    \toprule
    Model & Cut & OptimalTimes & GapWithoutInf & \#Inf \\
    \midrule
    DD-JS & IIS & 427.1 & 0.8 & 0 \\
    DD-JS & No-Good & 490.0 & 0.8 & 0 \\
    DD-LJ & IIS & 518.0 & 1.1 & 0 \\
    DD-LJ & No-Good & 547.3 & 1.0 & 0 \\
    IP & IIS & 136.7 & 1.6 & 101 \\
    IP & No-Good & 537.0 & 1.5 & 2 \\
    \bottomrule
    \end{tabular}
    \label{tab:ORSSummaryValues2}
	\end{center}
\end{table}

\begin{table}[ht]
	\begin{center}
	\caption{Detailed ORS Results}
    \begin{tabular}{llrrrr}
    \toprule
    Model & Cut & \#Callback & \#Cuts & ResolTime & ResolTimePerCB \\
    \midrule
    DD-JS & IIS & 112.4 & 18681.2 & 919.3 & 8.2 \\
    DD-JS & No-Good & 116.6 & 16235.7 & 925.4 & 7.9 \\
    DD-LJ & IIS & 88.0 & 14245.7 & 1005.4 & 11.4 \\
    DD-LJ & No-Good & 100.7 & 14025.8 & 976.9 & 9.7 \\
    IP & IIS & 4.6 & 189.7 & 48.2 & 10.5 \\
    IP & No-Good & 52.2 & 4190.3 & 779.3 & 14.9 \\
    \bottomrule
    \end{tabular}
    \begin{tabular}{llrrrr}
    \toprule
    Model & Cut & CreateCutTime & CreateSPTime & UpdateTime & ResolTime \\
    \midrule
    DD-JS & IIS & 0.5 & 0.4 & N/A & 919.3 \\
    DD-JS & No-Good & 0.5 & 0.3 & N/A & 925.4 \\
    DD-LJ & IIS & 0.4 & 1.9 & N/A & 1005.4 \\
    DD-LJ & No-Good & 0.4 & 1.8 & N/A & 976.9 \\
    IP & IIS & 909.0 & 247.4 & 13.1 & 48.2 \\
    IP & No-Good & 0.1 & 181.7 & 144.4 & 779.3 \\
    \bottomrule
    \end{tabular}
    \label{tab:ORSSummaryResultsDetails}
	\end{center}
\end{table}

\subsection{Equal} \label{Asub:EqualResults}

\begin{table}[H]
	\begin{center}
	\caption{Summary of Equal Unfiltered Results}
    \begin{tabular}{llrrrrr}
    \toprule
    Model & Cut & \#Problems & \#Optimal & TotalTime & Gap \\
    \midrule
    DD-JS & IIS & 135 & 84 & 709.9 & 0.4 \\
    DD-JS & No-Good & 135 & 72 & 758.8 & 0.6 \\
    DD-LJ & IIS & 135 & 67 & 825.2 & 0.6 \\
    DD-LJ & No-Good & 135 & 63 & 844.6 & 0.7 \\
    IP & IIS & 135 & 17 & 1121.1 & inf \\
    IP & No-Good & 135 & 49 & 957.6 & 1.1 \\
    \bottomrule
    \end{tabular}
	\label{tab:EqualSummaryValues}
	\end{center}
\end{table}

\begin{table}[H]
	\begin{center}
	\caption{Summary of Equal Filtered Results}
    \begin{tabular}{llrrrr}
    \toprule
    Model & Cut & OptimalTimes & GapWithoutInf & \#Inf \\
    \midrule
    DD-JS & IIS & 411.9 & 0.4 & 0 \\
    DD-JS & No-Good & 372.1 & 0.6 & 0 \\
    DD-LJ & IIS & 444.0 & 0.6 & 0 \\
    DD-LJ & No-Good & 415.3 & 0.7 & 0 \\
    IP & IIS & 222.9 & 2.7 & 58 \\
    IP & No-Good & 529.1 & 1.1 & 0 \\
    \bottomrule
    \end{tabular}
    \label{tab:EqualSummaryValues2}
	\end{center}
\end{table}

\begin{table}[ht]
	\begin{center}
	\caption{Detailed Equal Results}
    \begin{tabular}{llrrrr}
    \toprule
    Model & Cut & \#Callback & \#Cuts & ResolTime & ResolTimePerCB \\
    \midrule
    DD-JS & IIS & 76.9 & 14834.9 & 610.3 & 7.9 \\
    DD-JS & No-Good & 98.5 & 23944.3 & 648.7 & 6.6 \\
    DD-LJ & IIS & 63.7 & 12501.1 & 756.1 & 11.9 \\
    DD-LJ & No-Good & 88.1 & 21152.5 & 737.5 & 8.4 \\
    IP & IIS & 7.6 & 149.5 & 99.3 & 13.0 \\
    IP & No-Good & 43.1 & 5675.6 & 654.4 & 15.2 \\
    \bottomrule
    \end{tabular}
    \begin{tabular}{llrrrr}
    \toprule
    Model & Cut & CreateCutTime & CreateSPTime & UpdateTime & ResolTime \\
    \midrule
    DD-JS & IIS & 0.4 & 0.3 & N/A & 610.3 \\
    DD-JS & No-Good & 0.6 & 0.3 & N/A & 648.7 \\
    DD-LJ & IIS & 0.3 & 1.7 & N/A & 756.1 \\
    DD-LJ & No-Good & 0.6 & 1.6 & N/A & 737.5 \\
    IP & IIS & 719.3 & 251.9 & 26.9 & 99.3 \\
    IP & No-Good & 0.2 & 178.9 & 109.1 & 654.4 \\
    \bottomrule
    \end{tabular}
    \label{tab:EqualSummaryResultsDetails}
	\end{center}
\end{table}

\end{document}